\RequirePackage{etoolbox}
\csdef{input@path}{{style/}{graphics/}}
\documentclass[aap,MSNbibl,noautosecdot,dvips]{arximspdf}
\usepackage{url,breakurl}

\doi{10.1214/13-AAP1000} 
\volume{25}
\issue{2}
\pubyear{2015}
\firstpage{477}
\lastpage{522}

\makeatletter

\newcommand{\rrvert}{\vert}
\newcommand{\llvert}{\vert}
\def\cal{\mathcal}
\newcommand{\eqref}[1]{(\ref{#1})}
\def\xrightarrow{\rightarrow}

\def\Ge{\operatorname{Geom}}

\def\OL{\Omega(1)}
\def\a{\alpha}
\def\G{\Gamma}
\def\l{\lambda}
\def\p{\pi}
\def\s{\sigma}
\def\t{\tau}
\def\Om{\Omega}
\def\cB{{\cal B}}
\def\cT{{\cal T}}

\newtheorem{theorem}{Theorem}

\newtheorem{lemma}[theorem]{Lemma}

\newtheorem{corollary}[theorem]{Corollary}
\newproclaim{Remark}{Remark}
\newproclaim{definition}{Definition}
\def\cW{{\cal W}}
\def\cX{{\cal X}}
\newcommand{\ul}[1]{\mathbf{#1}}
\newcommand{\ull}[1]{\bolds{#1}}
\newcommand{\rdown}[1]{{\lfloor #1 \rfloor}}

\newcommand{\brac}[1]{(#1)}

\newcommand{\ra}{\rightarrow}
\newcommand{\rat}{\ra}

\newcommand{\set}[1]{\{#1\}}
\def\sm{\setminus}
\def\seq{\subseteq}

\def\E{\mathbf{E}}
\def\Pr{\operatorname{\mathbf{Pr}}}

\def\Upsi{\mathbf{I}}
\def\Row{\rho}

\newcommand{\ignore}[1]{}

\newcommand{\cA}{{\cal A}}

\def\OL{\Omega(1)}

\makeatother

\begin{document}
\begin{frontmatter}

\title{Viral processes by random walks on random regular graphs}
\runtitle{Viral processes on random regular graphs}

\begin{aug}
\author[A]{\fnms{Mohammed} \snm{Abdullah}\corref{}\thanksref{TT1}\ead[label=e1]{email.mohammed@gmail.com}},
\author[B]{\fnms{Colin} \snm{Cooper}\ead[label=e2]{colin.cooper@kcl.ac.uk}}
\and
\author[C]{\fnms{Moez} \snm{Draief}\thanksref{TT1}\ead[label=e3]{m.draief@imperial.ac.uk}}

\thankstext{TT1}{Supported by QNRF Grant NPRP 09-1150-2-448.}

\runauthor{M. Abdullah, C. Cooper and M. Draief}

\affiliation{University of Birmingham,
King's College London and
Imperial College London}

\address[A]{M. Abdullah\\
School of Mathematics\\
University of Birmingham\\
Edgbaston\\
Birmingham B15 2TT\\
United Kingdom\\
\printead{e1}}

\address[B]{C. Cooper\\
Department of Informatics\\
King's College London\\
Strand Building\\
London WC2R 2LS\\
United Kingdom\\
\printead{e2}}

\address[C]{M. Draief\\
Department of Electrical \\
\quad and Electronic Engineering\\
Imperial College London\\
South Kensington Campus\\
London SW7 2AZ\\
United Kingdom\\
\printead{e3}}
\end{aug}

\received{\smonth{7} \syear{2013}}
\revised{\smonth{10} \syear{2013}}

%
\begin{abstract}
We study the SIR epidemic model with infections carried by $k$
particles making independent random walks on a random regular graph.
Here we assume $k \leq n^\epsilon$, where $n$ is the number of vertices
in the random graph, and $\epsilon$ is some sufficiently small
constant. We give an edge-weighted graph reduction of the dynamics of
the process that allows us to apply standard results of Erd\H{o}s--R\'
enyi random graphs on the particle set. In particular, we show how the
parameters of the model give two thresholds: In the subcritical regime,
$O(\ln k)$ particles are infected. In the supercritical regime, for a
constant $\beta\in(0,1)$ determined by the parameters of the model,
$\beta k$ get infected with probability $\beta$, and $O(\ln k)$ get
infected with probability $(1-\beta)$. Finally, there is a regime in
which all $k$ particles are infected. Furthermore, the edge weights
give information about when a particle becomes infected. We exploit
this to give a completion time of the process for the SI case.
\end{abstract}

%


\end{frontmatter}

\section{\texorpdfstring{Introduction.}{Introduction}}\label{sec1}
The spread of an infection throughout a population, often referred to
loosely as an \emph{epidemic}, has come to be modeled in various ways
in the literature, spurred by the richness of domains in which the
notion of a virus has gone beyond the traditional biological
phenomenon. Electronic viruses over computer networks are not the only
extension; others include \emph{rumour spreading}~\cite{Pittel87} or
\emph{broadcasting} \cite{Frag} and \emph{viral marketing} \cite
{Amini}. Models may vary over domains, but the underlying principle is
one of spread of some unit of information or state through interaction
between individuals.

In much of the literature on the spread of epidemics as well as the
dissemination of information, individuals
reside at fixed vertices of a graph and the evolution of
the state of an individual depends on the state of its neighbours in
the graph. In particular if the graph is complete, mean-field
(homogeneous mixing) models have been exploited to study the outcome of
diffusion process \cite{DaleyGani01}. More recently, there has
been an increasing interest in understanding the impact of
the network topology on the spread of epidemics in networks with fixed
nodes; see
\cite{Draief-Massouli10} for a review of such results. There has,
however, been little analytical
work to date on models where the possible interactions between the
nodes are
dynamic; that is, the underlying network structure evolves in time.

We explore a particular instance of dynamic interaction by
assuming that individuals are mobile particles and
can only infect each other if they are in sufficiently close
proximity. The model is
motivated both by certain kinds of biological epidemics, whose
transmission may be
dominated by sites at which individuals gather in close proximity
(e.g., workplaces or
public transport for a disease like SARS, cattle markets for
foot-and-mouth disease,
etc.) and by malware. Furthermore, it is relevant to studying the
dissemination of information in opportunistic
networks~\cite{Chaintreaual07} where the information is transmitted
between users who happen to be in each other's range. As in the case of static
networks \cite{Pittel87}
one may be interested in the time it takes for the rumour to be
known to all users.

In our model (elaborated upon below) there are $k$ particles making
independent, discrete-time, synchronous random walks on an $n$-vertex
$r$-regular random graph~$G$. Each particle is in one of three states:
susceptible (S), infected (I), or recovered (R). An infected particle
can infect a susceptible particle, which remains infected for a fixed
\emph{infectious period} $\xi$ before recovering permanently. This is
known as the \emph{SIR epidemic model} and is extensively studied in
the literature. When $\xi=\infty$ (the \emph{SI model}) particles never
go from I to R.

Two questions can be asked: (1) When $\xi<\infty$, how many particles
ever get infected? (2) When $\xi=\infty$, what is the \emph{completion
time} of the process? That is, how long till the last infection? We
address both of these questions by reducing the dynamics of the process
to what we term an \emph{interaction graph}. This turns out to have the
structure of an Erd\H{o}s--R\'enyi (E--R) random graph $\mathcal
{G}_{k,\hat{q}}$ on the set of particles, where the edge probability
$\hat{q}$ is a function of the parameters of the model. Infected
particles are connected components in $\mathcal{G}_{k,\hat{q}}$, and so
well-known results from the literature on E--R random graphs can be
applied using our reduction to answer question (1). In particular, we
show how the parameters of the model produce two thresholds: In the
subcritical regime, $O(\ln k)$ particles are infected with high
probability (w.h.p.), that is, with probability tending to $1$ as $n
\rightarrow\infty$. In the supercritical regime, for a constant
$\beta
$ determined by the parameters of the model, $\beta k$ get infected
with probability $\beta$, and $O(\ln k)$ get infected with probability
$(1-\beta)$. Finally, there is a regime in which all $k$ particles are
infected w.h.p.

Furthermore, the interaction graph reduction assigns weights on the
edges that give information about when a particle becomes infected.
This information can be used for addressing question (2), which we do
by giving a convergence in probability. This is detailed in Theorem~\ref
{completiontimetheorem}.

While the metaphor of an epidemic is a motivating and illustrative one,
this work is
part of the more general scope of \emph{interacting particle systems};
see, for example, \cite{AlFi}, Chapter~14.


\section{\texorpdfstring{Model.}{Model}}
Let $r \geq3$ be fixed. Let $\mathcal{G}_{r}$ denote the set of $r$-regular
graphs with vertex set $V_G=\{1,2,\ldots,n\}$ and the uniform measure.
Let $G =(V_G,E_G)$ be chosen uniformly at random (u.a.r.) from
$\mathcal{G}_{r}$.
The results in this paper are
always asymptotic in $n=|V_G|$. The notation $O, o, \Omega, \Theta$
have their usual meanings with respect to $n$. We denote by $\Omega
^+(1)$ a quantity that can be replaced by an arbitrarily large
constant. In some contexts where the sign does not matter, we may write
$+o(1)$ instead of $-o(1)$, so as not to lend any significance to the
sign. It will be obvious when this is the case.

Denote the set of particles in the system by $\mathcal{P}$. At step
$t$, let $\mathcal{S}(t), \mathcal{I}(t), \mathcal{R}(t)$ be the set of
susceptible, infected, and recovered particles, respectively. Since a
given particle is in precisely one of these sets at a given time, they
form a partition of $\mathcal{P}$.

When two particles $x$ and $y$ meet at some vertex $v \in G$, we call
that an \emph{$xy$ meeting}, or say $x$ and $y$ are \emph{incident}.
When there is an $xy$ meeting at some time step $t$, an \emph
{interaction} takes place between them at that time step with
probability $\Row\in(0,1]$, which is a fixed constant parameter of
the model. We term such an event an \emph{$xy$ interaction} and call
$\Row$ the \emph{interaction probability}. If one of the particles is
infected and the other is susceptible, the infection takes place upon
interaction. The infectious period $\xi$ is not restricted to being
finite or constant; it is permitted to be $\infty$ or some function of
$n$, for example.

Consider that the time step counter has just changed from $t-1$ to $t$:
every particle $x$ makes an independent move in its random walk.
Subsequently, the rules are as follows:
\begin{longlist}[(1)]
\item[(1)] If $x$, $y$ are on the same vertex $v$, there is an $xy$
interaction with probability $\Row$; if they are on different vertices,
they cannot interact; each particle pair (non)interaction is
independent of every other particle pair (non)interaction.
\item[(2)] If $x \in\mathcal{S}(t-1)$, then $x \in\mathcal{S}(t)$,
unless there was an $xy$ interaction with at least one particle $y \in
\mathcal{I}(t-1)$. In the latter case, we say $x$ was infected at time
$t-1$ and write $t(x)=t-1$.
\item[(3)] If $x$ was infected at time $t(x)$, then $x \in\mathcal
{I}(t)$ for $t=t(x)+1, t(x)+2, \ldots, t(x)+\xi$. Subsequently, $x
\in
\mathcal{R}(t)$ for all $t > t(x)+\xi$.
\end{longlist}
[Note, we assume that $\mathcal{R}(0)=\varnothing$, and so a particle
could only be in $\mathcal{R}(t)$ if it had been infected at some time
previous to $t$.]

Observe two things from the above rules. Firstly, infections are not
transitive in a single time step. For example, suppose $x \in\mathcal
{I}(t-1)$ meets $y,z \in\mathcal{S}(t-1)$ at vertex $v$ at time step
$t$. If $x$ interacts with $y$ but not $z$, then $y$ does not pass on
the infection to $z$ at that time step, regardless of whether or not
they interact.

Secondly, the rules allow, in principle, for infections to be passed on
even with infectious period $\xi=1$; there could be a chain of
infections with $x$ infecting $y$, which in turn infects $z$, etc.

A note on notation and terminology: We write $[t_1, t_2]$ or $[t_1,
t_1+1, \ldots, t_2]$ to denote the set of time steps $\{t_1, t_1+1,
t_1+2 ,\ldots, t_2\}$, and we call these \emph{periods}. We may have
infinite periods, for example, $[T, T+1, \ldots]$. When we refer to
``time $t$'' we are referring to step $t$ on the counter---this is a
discrete-time process. The first time step is $t=1$, and $t=0$ is not a
time step, but a useful convention to express the initial condition.


\section{\texorpdfstring{Assumptions.}{Assumptions}}\label{assumptionssection}
We first specify some assumptions of the model.

If each particle is at distance at least $\omega (k,n)\equiv\Om( \ln
\ln n
+ \ln k)$ from every other, then we say the particles are in \emph
{general position} (g.p.).

We assume the following: {(i)} $G$ is typical (``typical'' is
defined in Section~\ref{TypicalGraphsSection}). {(ii)} The
number of particles $k \leq n^{\epsilon}$ where $\epsilon>0$ is some
sufficiently small constant. {(iii)} Particles start in general
position. {(iv)} $\mathcal{I}(0)=\{x_0\}$, and we refer to $x_0$
as the \emph{initial infective}.
Of course we also assume $|\mathcal{R}(0)|=0$.

A graph $G$ is \emph{typical} if it satisfies certain conditions. The
definition will be given in following sections, but for now it suffices
to say that most graphs in $\mathcal{G}_r$ are typical; that is, a
graph $G$ picked u.a.r. from $\mathcal{G}_r$ will be typical w.h.p.\looseness=-1

Assumption {(iii)} is also not unreasonable; it is
straightforward to verify that if the positions of each of the $k$
particles are chosen u.a.r. from the vertices of $G$, then w.h.p. they
will be in g.p. with respect to each other if $\epsilon$ is small enough.

It is not difficult to extend our results to a greater number of
initially infected particles, but we make assumption {(iv)} for
convenience and clarity.


\section{\texorpdfstring{Results.}{Results}}\label{resultssection}
Let $M_k$ be the total number of particles that ever get infected in
the course of the process, and let $T_k$ be the \emph{completion time}
of the process, the time step at which the last infection takes place.
Define
%
\begin{equation}
\theta_r = \frac{r-1}{r-2}\label{deftheta_r}
\end{equation}
and
%
\begin{equation}
\psi=\frac{\Row(r-1)}{r-2+\Row}\label{defpsi}
\end{equation}
(observe that $\Row\leq\psi\leq1$).

\begin{theorem}\label{epidemic}
Assume the conditions of Section~\ref{assumptionssection}, and suppose
$k \rightarrow\infty$ as $n \rightarrow\infty$. Let
%
\begin{equation}
\Phi= k \biggl(1- \biggl(1-\frac{\psi}{\theta_r n} \biggr)^\xi \biggr).
\label{np}
\end{equation}

\begin{longlist}[(iii)]
\item[(i)] If $\Phi<1$, then w.h.p., $M_k=O(\ln k)$.
\item[(ii)] If $\Phi\rightarrow c$ for any constant $c>1$, then with
probability $(1+o(1))\beta$, $M_k = (1+o(1))\beta k$, otherwise $M_k =
O(\ln k)$. $\beta$ is the unique solution in $(0,1)$ of the equation
$1-x=e^{-c x}$.
\item[(iii)] If $\Phi>(1+\varepsilon)\ln k$ where $\varepsilon>0$
is a
constant, then w.h.p. $M_k=k$ (i.e., all the particles get infected).
\end{longlist}
\end{theorem}

Theorem~\ref{epidemic} is for finite $\xi$, but observe that taking the
convention that $x^{\infty}=0$ for $|x|<1$ means that part
{(iii)} is consistent with the SI model, for which all particles get
infected almost surely. The theorem effectively gives conditions for
transitions between different possible ``regimes'' of behaviour. The
most interesting is the regime of part {(ii)}, which is entered
as $\Phi$ transitions from $\Phi<1$ to $\Phi>1$. Roughly speaking, in
this transition the number of infected particles goes from very few
($O(\ln k)$) w.h.p., to a constant fraction $\beta k$ with probability
$\beta$, or $O(\ln k)$ with probability $1-\beta$.

Concerning the completion times, we shall demonstrate that for $\xi=
\infty$, (i.e., the SI model) how the edge weightings can be exploited
by application of a theorem of \cite{Janson} to get the following.

\begin{theorem}\label{completiontimetheorem}
Assume the conditions of Section~\ref{assumptionssection}, and suppose
$k \rightarrow\infty$ as $n \rightarrow\infty$. When $\xi=\infty
$, we
have the following convergence in probability:
%
\begin{equation}
\frac{T_k}{n\ln k/k} \mathop{\xrightarrow}^{p} 2\frac{\theta_r}{\psi}, %
\end{equation}
where $T_k$ is the completion time for $k$ particles, that is, the time
at which the final particle is infected.
\end{theorem}


\section{\texorpdfstring{Related work.}{Related work}}\label{relatedwork}
In this section, we briefly describe some of the relevant related work
on diffusion processes like epidemic spreading and the dissemination
of information in mobile environments.

There has recently been a growing body of work in the interacting
particle systems
community analysing epidemic models with mobile particles. In
\cite{Dickman} the authors provide a review of the
results, techniques and conjectures when the graph is an infinite
lattice. In \cite{Rhodes}, the authors explore by means of
mean-field approximations the evolution of the number of infected
individuals when individuals perform random motions on the plane.
Recent papers by Peres et al. \cite{Peres},
Pettarin et al.~\cite{Upfal} and Lam et al. \cite{Lam} analyse mobile
networks modeled as multiple random walks; as Brownian motion on
$\mathbb{R}^d$ in \cite{Peres}, as walks on a $2$-dimensional grid in
\cite{Upfal} and as walks on a grids of dimension $3$ and above in
\cite
{Lam}. In each case, there is a parameter $r$ within which distance a
pair of walks can communicate, producing a communication graph (which
is disconnected below a certain threshold $r_c$). Peres et al. \cite
{Peres} study various aspects of the communication graph, such as how
long it takes a particular point to become part of the giant component.
Pettarin et al.~\cite{Upfal} study the broadcast time $T_B$ of a piece
of information and originating in one agent in relation to $r$. Setting
$r=0$ means information is passed only when particles are coincident.
In this case, $T_B$ is equivalent to our completion time, and the
authors of \cite{Upfal} give, for $k\geq2$, $T_B=\tilde{\Theta
}(n/\sqrt {k})$ w.h.p.\setcounter{footnote}{1}\footnote{The tilde notation hides polylogarithmic factors.
For example, $\tilde{O}(f(n))=O(f(n)\log^cn)$ for some constant $c$.}
In \cite{Lam} $r=1$ and the results they have show a significant
difference to the $2$-dimensional case of \cite{Upfal}. We state their
results for the $3$-dimensional case: There exists a constant $c$ such
that if $cn^{{1}/{3}}\log^2n^{{1}/{3}}<k<n$, then
$T_B=\tilde
{\Theta}(n^{{5}/{6}}/\sqrt{k})$; if $k<cn^{{1}/{3}}\log
^{-2}n^{{1}/{3}}$, then $T_B=\tilde{\Theta}(n/k)$. Statements are
w.h.p.

Of closer relevance to this work are \cite{Dimitriou} and \cite
{Draief}. In both of these papers, the authors study infections carried
by random walks on graphs. In particular, Dimitriou, Nikoletseas, and
Spirakis \cite{Dimitriou} analyse an SI model similar to ours; multiple
random walks on
a graph $G$ that carry a virus, and the quantity of interest is the
completion time. They give a general upper bound $\E[T_k]=O(m^* \ln k)$
for any graph $G$, where $m^*$ is the expected meeting time of a pair
of random walks maximised over all pairs of starting vertices. Special
cases are analysed too. In particular, they give an upper bound of $\E
[T_k]=O(\frac{nr}{k} \ln k \ln n)$ for random $r$-regular graphs. This
is a factor $\ln n$ larger than the precise value of the process
considered in this paper.

Finally, in \cite{Frag}, Baumann, Crescenzi, and Fraigniaud study
flooding on dynamic random networks. A fixed set of $n$ vertices is
part of a dynamic random graph process where each edge is an
independent two-state Markov chain, either existing or not existing. A
single initial vertex initially holds a message, and any vertex which
receives this message broadcasts it on existing edges for the next $k$ steps.
Although flooding is a different process to multiple random walks, the
authors develop a reduction to a weighted random graph with some
similarity to the interaction graphs we present. It allows them to
derive relations between the edge-Markov probabilities and state
asymptotic bounds on the number of vertices that receive the message,
as well as the broadcast (equivalently, completion) time.


\section{\texorpdfstring{Overview.}{Overview}}
Section~\ref{approachSection} addresses the behavior of (multiple)
random walks on a graph, drawing on established results from the
literature. It introduces the product graph framework, used for mapping
multiple walks on $G$ to a single walk on the product graph $H$. Thus,
analysis of the multiple walks can be done by analysis of the single
walk on $H$. After the background theory is given, the \emph{first
visit lemma} is presented. This lemma, first established in \cite
{CFreg}, then subsequently refined in other papers, gives the
probability of a walk visiting a vertex $v$ for the first time (after
mixing) at step $t$. It is used to establish new lemmas, created
specifically for the analysis of the problem in this paper. The
culmination is Lemma~\ref{ProbBxyLemma}, which calculates probabilities
of meetings of particle pairs. This lemma is the main tool of Section~\ref{approachSection} that is used in subsequent sections.

In Section~\ref{Interaction graph} we introduce the interaction graph
$\Upsi$, a complete grah on the particle set $\mathcal{P}$ with edge
weights determined by the outcome of the process. $\Upsi$ can therefore
be represented by a random ${k\choose2}$-vector, some entries of which
will be $\infty$ if their associated particles were never infected. The
edge weights give timing information on particle infection times,
encoding which particles get infected and when those infections occurred.

Our approach is to analyse special cases and build upon them until
reaching the full general case. We start in Section~\ref{SIModel} with
the SI model with $\Row=1$, then generalise that in Section~\ref{SIRRho1} to the SI(R), keeping $\Row=1$. Finally, in Section~\ref{TheGeneralCase}, we give the most general case, SI(R) with $\Row\leq
1$. We reiterate that whilst $\xi$ is permitted to be $\infty$ or take
any positive integer value (which may vary with $n$), $\Row$ is a constant.

A simple algorithm can be used to construct $\Upsi$ from the unfolding
of the process, but to actually calculate the probability of a
particular realisation of $\Upsi$, we will employ the tools of Section~\ref{approachSection}, in particular, Lemma~\ref{ProbBxyLemma}. This
lemma gives the probability of a first meeting time of a pair of
particles being at time $t$, but $t$ is required to be at least $\ell
=\Omega(T^3)$, where $T$ is the mixing time of the walks. It cannot
account for what happens in the first $\ell$ steps. As such, rather
than calculate the weights of $\Upsi$, we will couple the process to a
slightly modified version of it that demands that the interaction
probability $\Row$ is temporarily switched to zero in the ``blind''
periods---those length-$\ell$ periods which cannot be accounted
for---before being switched back. When $\Upsi$ is constructed under the
new process, it may turn out differently to what it would have been
under the original process. Therefore, for the $\Row=1$ case, we will
use $\Upsi'$ to specify the interaction graph under the modified process.

In Section~\ref{SIModel}, where $\xi=\infty$, it will be shown that the
edge weights of $\Upsi'$ (which will all be finite, almost surely) are
``almost'' independent, being well approximated by independent and
identically distributed (i.i.d.) random variables having geometric
distribution with parameter $1/(\theta_r n)$, denoted $\Ge(1/(\theta_r
n))$. Thus, $\Upsi'$ can be modeled by $\Lambda$, a complete graph with
i.i.d. random edge weights having that distribution. We can use
$\Lambda
$ to calculate distances and relate these to $\Upsi'$. We then need to
show $\Upsi'$ is a good approximation for $\Upsi$. We do so by showing
that w.h.p., the two graphs, constructed under their respective
processes, will give the same edge weights when the processes are coupled.

Section~\ref{SIModel} concludes with the proof of Theorem~\ref
{completiontimetheorem} for the special case $\Row=\psi=1$. This is
done in Section~\ref{comptimesection} by an application of a result of
Janson \cite{Janson} on distances in randomly edge-weighted graphs. The
theorem is applied to $\Lambda$ and transfers by the above arguments to
$\Upsi$.

In Section~\ref{SIRRho1}, where we address $\xi<\infty$ (but keep
$\Row
=1$), we try to quantify the size of the outbreak. A particle $x$ is
infected if and only if there is a $\xi$-path in $\Upsi$ from the
initial infective $x_0$ to $x$. A $\xi$-path is one in which all edge
weights in the path are at most $\xi$. Equivalently, the infected
particles form the connected component $x_0$ belongs to when we delete
all edges in $\Upsi$ with weights exceeding $\xi$. Referring to this
connected component as $\mathcal{C}_{f_\xi(\Upsi)}$, we relate
$\mathcal
{C}_{f_\xi(\Upsi)}$ to $\mathcal{C}_{f_\xi(\Lambda)}$, the equivalent
in~$\Lambda$. The function $f_\xi(F)$ takes a graph $F$ and returns the
same graph, but with edges weighing more than $\xi$ deleted. Ignoring
edge weights, $f_\xi(\Lambda)$ is an Erd\H{o}s--R\'enyi random graph on
$k$ vertices $\mathcal{G}_{k, \hat{q}}$ with edge probability $\hat
{q}=1-(1-\frac{1}{\theta_r n})^\xi$. Consequently, standard results on
Erd\H{o}s--R\'enyi random graphs give characterisations of the size of
$\mathcal{C}_{f_\xi(\Lambda)}$, which in turn transfer to $\mathcal
{C}_{f_\xi(\Upsi')}$ and subsequently to $\mathcal{C}_{f_\xi(\Upsi)}$.
Thus, we can determine how many particles get infected as a function of
the parameters of the model. This will give us Theorem~\ref{epidemic}
for the special case $\Row=1$.

We then move on to Section~\ref{TheGeneralCase}, where we deal with
$\Row<1$. We begin in Section~\ref{twoparticles} with a heuristic
treatment of a two particle system, consisting of the initial infective
$x_0$ and another particle, which is susceptible. This will allow us to
outline the techniques in a clear and concise way without being
hindered by detail. Subsequently, in Section~\ref{SIrholeq1case}, we
will formalise the arguments given in Section~\ref{twoparticles},
extended to all $k$ particles. The core of this section will be Lemma~\ref{mainPart3Lemma}. This lemma is essentially a generalisation of
Lemma~\ref{ProbBxyLemma}. It will allow us to determine probabilities
of having a first interaction at some time $t$, while allowing for the
possibility that there may have been meetings of particles prior to $t$
where no interaction took place. Section~\ref{SIrholeq1case} builds on
and generalises the previous ones where $\Row=1$, and we will detail
how Theorems \ref{epidemic} and~\ref{completiontimetheorem} are
justified in their full generality by the results in this section.

Finally, in Section~\ref{ConclusionSection}, we make concluding
remarks, including possible extensions.

\section{\texorpdfstring{Random walks on graphs: Tools and techniques.}
{Random walks on graphs: Tools and techniques}}\label{approachSection}
In this section we detail key concepts and lemmas that we use to
analyse the viral process. In Section~\ref{TypicalGraphsSection} we
give a formal definition of typical graphs. Knowing the properties of
typical graphs, we can make statements about how walks behave on them.
In Section~\ref{Convergence to stationarity and product graph
formulation section} we describe how the long-term behaviour of a
random walk---specifically its convergence to the stationary
distribution---relates to the eigenvalues of its transition matrix. We
also show how to map the $k$ multiple walks on $G$ on to a single walk
on another graph, the \emph{product graph} $H$. Much of the analysis of
meeting times between walks is done through the framework of the
product graph. A particular pair of particles meeting in $G$ maps to
the single walk on $H$ being at a set of vertices of $H$. Such a set of
vertices is contracted to a single vertex, resulting in a derived graph
$\Gamma$ upon which a single walk moves. Thus calculations of particle
meeting times in $G$ are done by calculating the hitting times of a
single walk to a contracted vertex in $\G$.

Lemma~\ref{T-value} shows that for each of the graphs $G, H, \G$, the
walks are \emph{rapidly mixing}, meaning that they converge to their
respective stationary distributions quickly. This is a crucial
component of the proofs and the behaviour of the processes. In Section~\ref{ssec6} we introduce a key lemma, Lemma~\ref{L3}, which allows us
to make precise calculations of $ f_t(u \rat v)$, the probability that
a walk starting at $u$ visits $v$ for the first time (after mixing) at
time $t$. When this lemma is applied to contracted vertices in $\G$, it
gives us probabilities of meeting times of particles walking on $G$.
Lemma~\ref{contract} in Section~\ref{ssec6} gives formal justification
for the use of vertex contraction to reason about visits to sets, and
Section~\ref{Which vertex section} investigates the probability that a
particular vertex $v$ was visited when a set of vertices $S$ with $v
\in S$ has been visited.

Lemma~\ref{L3}, in conjunction with Lemmas \ref{contract}, \ref
{whichvx}, and \ref{pALem} culminate in Lemma~\ref{ProbBxyLemma} in
Section~\ref{Particle pair meetings section}. This lemma gives
probabilities for a particular pair of particles meeting at time $t$
and no other pair meeting before hand. It is the main tool used in
Section~\ref{SIModel}.

\subsection{\texorpdfstring{Typical graphs.}{Typical graphs}}\label{TypicalGraphsSection}
Let
%
\begin{equation}
\label{sig} L_1= \lfloor \epsilon _1
\log_rn\rfloor,
\end{equation}
where $\epsilon _1>0$ is a sufficiently small constant.

A vertex $v$ is \emph{treelike} if there
is no cycle in the subgraph $G[v,L_1]$ induced by the set of
vertices within (graph) distance $L_1$ of $v$.

A cycle $C$ is \emph{small} if $|C|\leq L_1$.

\begin{longlist}[P4.]
\item[P1.] $G$ is connected and not bipartite.
\item[P2.] The second eigenvalue of the adjacency matrix of $G$
is at most
$2\sqrt{r-1}+\varepsilon$, where $\varepsilon>0$ is an arbitrarily
small constant.
\item[P3.] There are at most $n^{2\epsilon _1}$ vertices on
small cycles.
\item[P4.] No pair of cycles $C_1,C_2$ with $|C_1|,|C_2|\leq100L_1$
are within distance $100L_1$ of each other.
\end{longlist}

We say an $r$-regular graph $G$ is \emph{typical} if it satisfies
properties P1--P4.

Note that P3 implies that at most $n^{\epsilon _C}$
vertices of a typical $r$-regular graph are not treelike, where
%
\begin{equation}
\label{nontree} n^{\epsilon _C}=O \bigl(r^{L_1} n^{2\epsilon _1}
\bigr)=O \bigl(n^{3
\epsilon _1} \bigr).
\end{equation}

\begin{lemma}[(\cite{CFR-Mult})]\label{typG}
Let ${\cal G}_r' \seq{\cal G}_r$ be the set of typical $r$-regular
graphs. Then
$|{\cal G}_r'|/|{\cal G}_r| \rightarrow1 $ as $n \rightarrow\infty$.
\end{lemma}

{P1} implies that a random walk will converge to a stationary
distribution $\pi$ on the vertex set. Because the graph is regular,
$\pi
$ will be the uniform distribution. {P2} implies that a random walk
will converge quickly to the stationary distribution; it will be \emph
{rapidly mixing}. In fact, for all the graphs we consider, $O(k \ln n)$
steps is sufficient for our results. A typical graph also has mostly
treelike vertices.

\subsection{\texorpdfstring{Convergence to stationarity and product graph
formulation.}{Convergence to stationarity and product graph
formulation}}\label{Convergence to stationarity and product graph
formulation section}

Let $G$ be a connected graph with $n$ vertices and $m$ edges.
For random walk $\cW_{u}$ starting at a vertex $u$ of $G$,
let $\cW_{u}(t)$ be the vertex
reached at step $t$. Let $P=P(G)$ be the matrix of
transition probabilities of the walk, and let
$P_{u}^{t}(v)=\Pr(\cW_{u}(t)=v)$.
If the random walk $\cW_{u}$ on $G$ is ergodic, it will converge to
stationary distribution~$\pi$.
Here $\pi(v)=d(v)/(2m)$, where $d(v)$ is the degree of vertex $v$.
We often write $\pi(v)$ as $\pi_v$.
The eigenvalues of $P(G)$ are $\l_0=1\ge
\l_1\ge\cdots\ge\l_{n-1} \geq-1$. Let
$\l_{\max}=\max(\l_1,|\l_{n-1}|)$.
The rate of convergence of the walk is given by
%
\begin{equation}
\label{mix}\bigl |P_{u}^{t}(x)-\pi_x \bigr| \leq(
\p_x/\p_u)^{1/2}\l_{\max}^t.
\end{equation}
For a proof of this, see, for example, Lov\'{a}sz \cite{Lo}.

To ensure that the walk is both ergodic and that $\l_{\max}=\l_1$, we
make the chain
\emph{lazy}; that is, the walk only moves to a neighbour with
probability $1/2$.
Otherwise it stays where it is. This shifts each eigenvalue up by $1$,
and so $\lambda_1 \geq\lambda_{n-1}=0$,
and~\eqref{mix} still holds.

Now define the \emph{product graph} $H=H(G,k)=(V_H,E_H)$ to
have vertex set $V_H=V^k$
and edge set $E_H=E^k$.
The vertices $\ul v$ of $H$ consist of $k$-tuples
$\ul v= (v_1,v_2,\ldots,v_k)$ of vertices $v_i \in V_G, i=1,\ldots,k$,
with repeats allowed.
Two vertices $\ul v, \ul w$ are adjacent if $(v_1,w_1),
\ldots,(v_k,w_k)$ are edges of $G$. The purpose of defining the graph
$H$ is that
we can replace the $k$ random walks $\cW_{u_i}(t)$ on $G$ with current
positions $v_i$
and starting positions $u_i$ by a single walk $\cW^H_{\ul u}(t)$. Note
that because $G$ is assumed to be simple, no vertex in $V_G$ has a
loop. Consequently, no edge $e=((v_1,w_1),
\ldots,(v_k,w_k))$ in
$E_H$ has $v_i=w_i$ for any $1 \leq i \leq k$. This is the case despite
the actual walk on the graph being ``lazy'' for part of the time, as
will be described below.

We introduce some extra notation: for a graph $F=(V_F, E_F)$, a vertex
$v \in V_F$, and a set of vertices $S \subseteq V_F$, let $d_F(v)$ be
the degree of vertex $v$ in $F$, and let $d_F(S)=\sum_{v \in S}d_F(v)$.

Now, if $S\seq V_H$, then $\Gamma=\G(S)$ is obtained from $H$ by contracting
$S$ to a single vertex $\gamma (S)$. All edges, including loops and parallel
edges are retained, producing a multigraph.
Thus $d_{\G}(\gamma )=d_H(S)=r^k|S|$.
Moreover $\G$ and $H$ have the same total degree $(nr)^k$, and
the degree of any vertex of $\G$, except $\gamma $, is $r^k$.

\begin{definition}[({Mixing time, maximal mixing time})]
For $F\in\{G,H,\G\}$, let $\cW^F_{u}$
be a lazy random walk starting at $u \in V_F$. The \emph{mixing time}
$T_F$ is the smallest $t$ such that, for graph
$F=(V_F,E_F)$ and $t\geq T_F$, 
\begin{equation}
\max_{u,x\in V_F} \bigl\llvert P_{u}^{t}(x)-
\pi_x \bigr\rrvert \leq\label{newDistance} \frac{\min_{x \in V_F}(\pi_x)}{n^3}.
\end{equation}
A \emph{maximal mixing time} $T$ is defined as
\[
T = \max \bigl\{T_{\G(S)} \dvtx S\seq V_H d_H(S)
\le k^2 n^{k-1}r^k \bigr\}.
\]
\end{definition}

Observe that the mixing time $T_\Gamma$ depends on the particular set
of vertices $S$ that gets contracted. The maximal mixing time is
defined for convenience; if there is ambiguity about what the mixing
time is, or a single mixing time is stated for a number of contractions
of sets $S$, then it is safe to assume the maximal mixing time.

The following is a slightly modified version of a lemma proved in \cite
{CFR-Mult}.

\begin{lemma}[(\cite{CFR-Mult})]\label{T-value}
Let $G$ be typical, and let $S\seq V_H$ be such that $d_H(S)\le k^2 n^{k-1}r^k$.

For $k \le n$,
\[
T_G=O(\ln n),\qquad T_H=O(\ln n)\quad \mbox{and}\quad
T_{\G}=O(k\ln n).
\]
\end{lemma}

As a consequence, a maximal mixing time $T$ has $T=O(k \ln n)$.

We analyse our walks in the product graph and assume that
we keep the chain lazy for the duration of the mixing time. At this
point it is mixed,
and we can stop being lazy. A lazy walk in the product graph maps to a
process where all the walks move or do not move together.
That is, with probability $1/2$, each walk independently takes a random
step, and with probability $1/2$ none of them do.
Consider the following two conditions: (i) interactions can only take
place upon moving to a new vertex, and (ii) $\xi$ can only be decreased
(by $1$) upon moving to a new vertex. It is not difficult to see (e.g.,
through coupling) that the laziness of the walk cannot affect the
infection outcomes. Laziness affects time, but only during mixing
periods since we do not keep the chain lazy thereafter.

The following lemma formalises the notion that we can deal with a first
visit (after the mixing time) to a member of a set $S$ of vertices of a
graph $H$ by contracting $S$ into a single vertex $\gamma=\gamma(S)$
and instead deal with a first visit to $\gamma$ on this altered graph.

\begin{lemma}[(\cite{CFgiant})]\label{contract}
Let $G$ be typical, $S\seq V_H$ be such that $d_H(S)\le k^2
n^{k-1}r^k$, and let $k \leq n^\epsilon $ for sufficiently small
$\epsilon $.

Let $\cW^H_u$ be a random walk in $H$ starting at $u \notin S$, and
let $\cW^{\Gamma}_u$ be a random walk in
$\G$ starting at the same vertex $u \ne\gamma $.
Let $T$ be a mixing time satisfying \eqref{newDistance} in both $H$ and
$ \G$. Let $\cA_w(t)$ be the event that
no visit was made to $w$ in the period $[T,T+1,\ldots,t]$. Then
\[
\Pr \bigl(\cA_{\gamma }(t);\G \bigr)=\Pr \biggl(\bigwedge_{v \in S}
\cA_{v}(t);H \biggr) \biggl(1 +O\biggl(\frac{1} {n^3}\biggr)
\biggr),
\]
where the probabilities are those derived from the walk in the given graph.
\end{lemma}

\subsection{\texorpdfstring{First visit lemma.}{First visit lemma}}\label{ssec6}
In this section, we introduce existing results about ``first visit''
behaviours for random walks on graphs. For the graphs we consider,
Lemma~\ref{L3} below, the \emph{first visit lemma} was proved in
\cite
{CFR-Mult}. It was initially presented in \cite{CFreg} then refined and
applied in a series of subsequent papers, amongst them \cite{CFgiant,CFR-Mult}.

Below, for some $\epsilon>0$, we denote by the term $\Omega^+(1)$ a
quantity that is at least some positive constant $C(\epsilon)$, which
can be made arbitrarily large by making $\epsilon$ sufficiently small.

\begin{lemma}[(First visit lemma \cite{CFreg,CFR-Mult})]\label{L3}
Let $G$ be typical, $S\seq V_H$ be such that $d_H(S)\le k^2
n^{k-1}r^k$, and let $k \leq n^\epsilon $ for sufficiently small
$\epsilon $.

For $F\in\{G,H,\G\}$ let $T=T_F$ and $\cW^F_u$ be a walk started at $u
\in V_F$. Let $f_t=f_t(u \rat v)$
be the probability that the first visit of $\cW^F_u$ to $v \in V_F$ in
the period $[T,T+1,\ldots]$ occurs at step $t$, and let $\cA_{v}(t)$ be
the event that $\cW_u$ does not visit $v$ in the period $[T,T+1,\ldots,t]$.

Let
%
\begin{equation}
\lambda=\frac{1}{KT} \label{lamby}
\end{equation}
for some sufficiently large constant $K$.
Then,
for all $t\geq T$,
%
\begin{equation}
\label{frat} f_t(u \rat v)= \bigl(1+O(T\pi_v) \bigr)
\frac{p_v}{(1+p_v)^{t+1}}+ O \bigl(T\pi_ve^{-\l t/2} \bigr)
\end{equation}
and
%
\begin{equation}
\label{atv} \Pr \bigl(\cA_{v}(t) \bigr)= \frac{(1+O(T\p_v))} {
{(1+(1+O(T\p_v))\p_v/R_v)}^t} +O
\bigl(T^2\pi_ve^{-
\lambda t/2} \bigr),
\end{equation}
where
%
\begin{equation}
\label{pv} p_v=\frac{\pi_v}{R_v(1+O(T\pi_v))}
\end{equation}
and
%
\begin{equation}
\label{Rtheta} R_v = \theta_r+O
\bigl(k^2n^{-\Omega(1)}+(k \ln n)^{-\Omega^+(1)} \bigr) .
\end{equation}
\end{lemma}

We briefly discuss the terms $R_v$ in \eqref{pv} and $\Omega^+(1)$. For
a given graph, $R_v$ is the expected number of returns in the mixing
time to a vertex $v$, for a walk that starts at $v$. The initial
placement of the walk at $v$ at time $t=0$ is counted. To usefully
apply the lemma, one needs to calculate (or approximate) $R_v$.
Consider the case where $v$ is a vertex of a $r$-regular random graph
$G$. The walk is rapidly mixing [the mixing time $T$ being $O(\ln n)$],
and for most vertices $v$, the local structure is a tree. Because of
this, the quantity $\theta_r=\frac{r-1}{r-2}$, the expected number of
returns (ever) to the root of an $r$-regular infinite tree, provides a
close approximation to $R_v$ on~$G$. For $R_{\gamma(S)}$, the expected
number of returns in the mixing time to a contracted vertex $\gamma(S)$
in the product graph, our bound on the mixing time $T_\Gamma$ is $O(k
\ln n)$. However, it turns out that $\theta_r$ also provides a good
approximation for $R_{\gamma(S)}$. This was determined in \cite
{CFR-Mult}, and our statement of the first visit lemma incorporates
both the general statement introduced in \cite{CFreg}, as well as the
bound on $R_{\gamma(S)}$ given in \cite{CFR-Mult}. The fact that the
$\Omega^+(1)$ term in \eqref{Rtheta} can be an arbitrarily large
constant is demonstrated in the derivation of $R_v$ in \cite{CFR-Mult},
Lemma~19.

A consequence of the condition $d_H(S)\le k^2 n^{k-1}r^k$ is that
$d(\gamma(S))\leq\break k^2 n^{k-1}r^k$, so $\pi_{\gamma(S)}\leq k^2/n$.
Since $R_{\gamma(S)}= (1+o(1))\theta_r$, we have $p_{\gamma
(S)}=O(\pi
_{\gamma(S)}/R_{\gamma(S)})=O(\pi_{\gamma(S)})=O(k^2/n)$. Therefore,
$T\pi_{\gamma(S)}=O(k^3 \ln n/n)=o(1)$, if $k \leq n^\epsilon$ and
$\epsilon$ is small enough.

We will rewrite \eqref{frat} and \eqref{atv} in a form that is more
natural in the context of this paper, in particular, a form that
resembles that of a geometric distribution. First note that
\[
\bigl((1-p_v) (1+p_v) \bigr)^t=
\bigl(1-p_v^2 \bigr)^t=1-O
\bigl(p_v^2t \bigr),
\]
so
\begin{eqnarray*}
\frac{1}{(1+p_v)^t}&=&\frac{(1-p_v)^t}{1-O(p_v^2t)}
\\
&=& \bigl(1+O \bigl(p_v^2t \bigr) \bigr)
(1-p_v)^t
\\
&=& \bigl(1+O \bigl(\pi_v^2t \bigr) \bigr)
(1-p_v)^t.
\end{eqnarray*}
In the above we have used the fact that $p_v=\Theta(\pi_v)$ since $R_v
= (1+o(1))\theta_r$ for the graphs in this paper.

We re-write \ref{frat} as follows:
\begin{eqnarray*}
f_t(u \rat v)&=& \bigl(1+O(T\pi_v) \bigr)
\frac{p_v}{(1+p_v)^{t+1}}+O \bigl(T\pi _ve^{-\l t/2} \bigr)
\\
&=& \bigl(1+O(T\pi_v) \bigr) \bigl(1-O(p_v) \bigr)
\bigl(1+O \bigl(\pi_v^2t \bigr) \bigr)p_v(1-p_v)^t
\\
&&{}+ O \bigl(T\pi_ve^{-\l t/2} \bigr)
\\
&=& \bigl[1+O(T\pi_v)+O \bigl(\pi_v^2t
\bigr)+O \bigl(Te^{-\l
t/2}(1-p_v)^{-t} \bigr)
\bigr]p_v(1-p_v)^t.
\end{eqnarray*}

Now
\begin{eqnarray*}
e^{-\lambda t/2}(1-p_v)^{-t}&=& \bigl(e^{{1}/{(2KT)}}(1-p_v)
\bigr)^{-t}
\\
&\leq& \bigl(e^{{1}/{(2KT)}}e^{-2p_v} \bigr)^{-t},
\end{eqnarray*}
and $\frac{1}{2KT}-2p_v=\Omega(1/T)$ since, as discussed above, $T\pi
_v=o(1)$.

Thus, we can write \eqref{frat} in the form
%
\begin{equation}
\label{frat2} \qquad f_t(u \rat v)= \bigl[1+O(T\pi_v)+O
\bigl(k^4t/n^2 \bigr)+O \bigl(T^2
e^{-\Omega
(t/T)} \bigr) \bigr]p_v(1-p_v)^t.
\end{equation}

Similarly, \eqref{atv} can be written as
%
\begin{equation}
\label{atv2} \Pr \bigl(\cA_{v}(t) \bigr)= \bigl[1+O(T
\pi_v)+O \bigl(k^4t/n^2 \bigr)+O
\bigl(T^2 e^{-\Omega
(t/T)} \bigr) \bigr](1-p_v)^t.
\end{equation}

\subsection{Which vertex in the set $S$ was visited?}\label{Which
vertex section}
For a set of vertices $S$, the following lemma gives the probability
that a particular vertex $v \in S$ is visited when $S$ is visited for
the first time after the mixing time.

\begin{lemma}\label{whichvx}
Let $G$ be typical, $S\seq V_H$ be such that $d_H(S)\le k^2 n^{k-1}r^k$
and let $k \leq n^\epsilon $ for sufficiently small $\epsilon $.

Let $\gamma $ be the contraction of $S$ in $H$, and for $v \in S$, let
$\delta$ be the contraction of $S\sm\{v\}$ in $H$, resulting in graphs
$\Gamma(S)$ and $\G(S\sm\{v\})$, respectively.
Let $T$ be the mixing time satisfying \eqref{newDistance} in both
$\Gamma(S)$ and $\Gamma(S \setminus\{v\})$.
Let $p_v, p_\gamma$ and $p_\delta$ be as given by \eqref{pv} for
$v,\gamma,\delta$ in their respective graphs.
Let $\epsilon _v$ be the solution to $p_{\gamma }-p_{\delta}=p_v
(1+\epsilon _v)$.

For $t\ge2(T+L)$ where $L=T^3$, let $ \cB_v=\cB_v(t)$ be the event
that the first visit to $S$ in the period $[T, T+1, \ldots]$ occurs at
step $t$ and that the visit is to node $v\in S$. Then
%
\begin{equation}
\Pr(\cB_v )= \biggl(1+ \bigl(1+o(1) \bigr)\epsilon
_v+O(L\pi_{\gamma
})+O \biggl(\frac
{k^4t}{n^2} \biggr)
\biggr)p_v (1-p_{\gamma })^t. \label{probB}
\end{equation}
\end{lemma}

Note that when $v$ is connected only to other vertices in $S$, it must
be that $\Pr(\cB_v )=0$. The RHS of \eqref{probB} is consistent with
this since in such a case, $\epsilon _v = -1\pm o(1)$.

\begin{pf*}{Proof of Lemma~\ref{whichvx}}
It is enough to prove the lemma for a two-vertex set $S=\{u,v\}$, as
one vertex
can always be a contraction of a set.
Let $t$ be expressed as $t=2T+L+s$, where $s \ge L$.
Divide $[0,t]$ into successive intervals of length $T,s,T,L$,
respectively, that is,
$[0,T-1], [T,s+T-1],
[s+T,s+2T-1], [s+2T,t]$.

Let $\cA$ be the event that $\cW(\s) \notin\{u,v\}$ for $ \s\in
[T,s+T-1]$
and that $\cW(t)=u$, but $\cW(\s) \ne u$ for $ \s\in[s+2T,t-1]$.
Contract $S$ to make $\gamma =\gamma (S)$ in $[T,T+s-1]$. Applying
\eqref{atv2}
to the period $[T,s+T-1]$, that is, letting $t=s+T-1$, and noting that
$T=O(k \ln n)$ and $s \geq L = T^3$
\begin{eqnarray*}
&&\Pr \bigl(\cA_{\gamma}(s+T-1) \bigr)\\
&&\qquad= \biggl(1+O(T\p_\gamma
)+O \biggl(\frac
{k^4s}{n^2} \biggr)+O \bigl(T^2 e^{-\Omega (T^2 )}
\bigr) \biggr) (1-p_\gamma)^{s+T-1}
\\
&&\qquad= \biggl(1+O(T\p_\gamma)+O \biggl(\frac{k^4s}{n^2} \biggr) \biggr)
(1-p_\gamma)^{s+T-1}.
\end{eqnarray*}

Now, starting from some vertex $x$ at time $s+T$, we apply \eqref
{frat2} for $u$ to not be visited in the period $[s+2T, t-1]$ then be
visited at $t$,
\begin{eqnarray*}
&&f_{t-s-T}(x \rat u)\\
&&\qquad= \biggl(1+O(T\p_u)+O \biggl(
\frac{k^4t}{n^2} \biggr)+O \bigl(T^2 e^{-\Omega ((t-s-T)/T )} \bigr) \biggr)
(1-p_u)^{t-s-T}p_u.
\end{eqnarray*}
Noting $t-s-T=L+T\geq T^3$,
\[
f_{t-s-T}(x \rat u)= \biggl(1+O(T\p_u)+O \biggl(
\frac{k^4t}{n^2} \biggr) \biggr) (1-p_u)^{L+T}p_u.
\]

Multiplying them together,
\begin{eqnarray*}
\Pr(\mathcal{A})&=&\Pr \bigl(\cA_{\gamma}(s+T-1) \bigr)f_{t-s-T}(x
\rat u)
\\
&=& \biggl(1+O(T\p_\gamma)+O \biggl(\frac{k^4t}{n^2} \biggr) \biggr)
(1-p_\gamma)^{s+T-1}(1-p_u)^{L+T}p_u
\\
&\leq& \biggl(1+O(T\p_\gamma)+O \biggl(\frac{k^4t}{n^2} \biggr)
\biggr) (1-p_\gamma)^{s}(1-p_u)^{L}p_u.
\end{eqnarray*}

Now,
\begin{eqnarray*}
(1-p_\gamma)^{s}(1-p_u)^L&=&(1-p_\gamma)^t
\biggl(\frac
{1-p_u}{1-p_\gamma
} \biggr)^L(1-p_\gamma)^{-2T}
\\
&=&(1-p_\gamma)^t \biggl(1+\frac{p_\gamma-p_u}{1-p_\gamma}
\biggr)^L \bigl(1+O(T\p_\gamma) \bigr)
\\
&=&(1-p_\gamma)^t \bigl(1+O(L\p_\gamma) \bigr)
\bigl(1+O(T\p_\gamma) \bigr)
\\
&=&(1-p_\gamma)^t \bigl(1+O(L\p_\gamma) \bigr).
\end{eqnarray*}

Let $\cB_u$ be the event that $\cW(t)=u$ and
$\cW(\s) \notin\{u,v\}$ for $\s\in[T,t-1]$.
Then $\cB_u \seq\cA$ and so $\Pr(\cB_u) \le\Pr(\cA) $.
It follows that
%
\begin{equation}
\label{Buu} \Pr(\cB_u) \le p_u(1-p_{\gamma })^t
\biggl(1+O(L\pi_{\gamma
})+O \biggl(\frac
{k^4t}{n^2} \biggr) \biggr).
\end{equation}

However, by contracting $S$ we have that
\[
\Pr(\cB_u \cup\cB_v)= \bigl(1+O(T
\pi_{\gamma })+O \bigl(k^4t/n^2 \bigr) \bigr)
p_{\gamma } (1-p_{\gamma })^t,
\]
and so
%
\begin{eqnarray}\label{Bvv}
\Pr(\cB_v) & \ge&\Pr(\cB_u \cup\cB_v)-
\Pr(\cB_u)
\nonumber
\\
&\ge& \bigl(1+O(L\pi_{\gamma })+O \bigl(k^4t/n^2
\bigr) \bigr) (p_{\gamma
}-p_u) (1-p_{\gamma })^t
\\
&=& \bigl(1+ \bigl(1+o(1) \bigr)\epsilon _v+O(L
\pi_{\gamma
})+O \bigl(k^4t/n^2 \bigr)
\bigr)p_v(1-p_{\gamma })^t.\nonumber
\end{eqnarray}
The result follows from \eqref{Buu} and \eqref{Bvv}.
\end{pf*}

\subsection{\texorpdfstring{Particle pair meetings.}{Particle pair meetings}}\label{Particle pair meetings section}
Consider the (unordered) pair of particles $(x,y)$, $x,y=1,\dots,k$.
Particles $x$ and $y$ being at the same vertex in $G$ maps in the
product graph to a set of vertices $S=\{ \ul v= (v_1,v_2,\ldots,v_k) \dvtx v_x=v_y\}\subset V_H$. We can, therefore, calculate the probability of
an $xy$ meeting in $G$ at time $t$ by calculating the probability of
the single random walk $\cW^H_{\ul u}$ on $H$ visiting $S$ at time $t$.
This in turn is done by the contraction described above, and
calculating the probability of the walk $\cW^{\Gamma}_{\ul u}$ on
$\Gamma$ visiting $\gamma (S)$ at time $t$. These two times are
asymptotically equal by Lemma~\ref{contract}.

More generally, let $A \subseteq\{(x,y)\dvtx x,y \in\mathcal{P}, x \neq
y\}$ be a set of particle pairs.
Consider the event $\{\mbox{for some $(x,y) \in A$ there is an $xy$
meeting at time $t$}\}$. In the product graph this maps to the event $\{
\cW^H_{\ul u}(t) \in S\}$ where $S=\{ \ul v= (v_1,v_2,\ldots,v_k) \dvtx v_x=v_y \mbox{ for some } (x,y) \in A\}$.

To use Lemma~\ref{L3} for the walk $\cW^{\Gamma}_{\ul u}$ on $\Gamma$
visiting $\gamma (S)$, we need
to calculate the relevant $p_v$ as per \eqref{pv}.

\begin{lemma}\label{pALem}
Let $G$ be typical and let $k \leq n^\epsilon $ for sufficiently small
$\epsilon $.

For a set of particle pairs $A$, let $S(A) \subset V_H$ be such that
${\ul v}=(v_1, v_2, \ldots,\break v_k) \in S$ if and only if, for some pair
$(x,y) \in A$, $v_x=v_y$ where $v_x$ (resp., $v_y$) is the position of
particle $x$ (resp., $y$) in $G$.
Then in $\Gamma=\Gamma(S(A))$,
%
\begin{equation}
p_{\gamma}=\frac{|A|}{\theta_r n} \biggl(1-O \biggl(\frac
{1}{n^{\Omega
(1)}}+
\frac{1}{(k \ln n)^{\Omega^+(1)}} \biggr) \biggr),
\end{equation}
where $\gamma=\gamma(S(A))$.
\end{lemma}

\begin{pf}
Let $N=|S(A)|$. A particular pair $(x,y) \in A$ can be on $n$ possible
different vertices in $G$, and for each one, the other particles can be
on $n^{k-2}$. Thus $N \leq|A|n^{k-1}$. Further, $N\geq N'$ where $N'$
is the number of $\ul v \in H$ such that only one of the particle pairs
occupy the same node of the graph, and
\[
N' \geq|A|n^{k-1}- |A|^2n^{k-2}=
|A|n^{k-1} \biggl(1-O \biggl(\frac
{k^2}{n} \biggr) \biggr).
\]

Thus, $N = |A|n^{k-1} (1-O (\frac{k^2}{n} ) )$ and
since each vertex of $H$ has degree $r^k$ and contraction preserves degree,
\[
\pi_{\gamma} = \frac{|A|r^kn^{k-1}}{n^kr^k} \biggl(1-O \biggl(\frac
{k^2}{n}
\biggr) \biggr) = \frac{|A|}{n} \biggl(1-O \biggl(\frac
{k^2}{n} \biggr)
\biggr)
\]
and
\[
T\pi_{\gamma} =O \biggl(\frac{|A|k \ln n}{n} \biggr) \biggl(1-O \biggl(
\frac
{k^2}{n} \biggr) \biggr) = O \biggl(\frac{k^3 \ln n}{n} \biggr)
\]
since $|A|\leq{k\choose2}$.\vspace*{1pt}

Hence \eqref{pv} becomes
\begin{eqnarray*}
p_{\gamma}&=& \frac{({|A|}/{n}) (1-O ({k^2}/{n}
) )}{(\theta_r+O(k^2n^{-\Omega(1)} +(k \ln n)^{-\Omega
^+(1)})
(1+O ({k^3 \ln n}/{n} ) )}
\\
&=&\frac{|A|}{\theta_r n} \biggl(1-O \biggl(\frac{k^3 \ln n}{n^{\Omega
(1)}}+\frac{1}{(k \ln n)^{\Omega^+(1)}}
\biggr) \biggr)
\\
&=&\frac{|A|}{\theta_r n} \biggl(1-O \biggl(\frac{1}{n^{\Omega
(1)}}+\frac
{1}{(k \ln n)^{\Omega^+(1)}}
\biggr) \biggr)
\end{eqnarray*}
when $\epsilon$ is small enough.
\end{pf}

\begin{lemma}\label{ProbBxyLemma}
Let $G$ be typical, and let $k \leq n^\epsilon $ for sufficiently
small $\epsilon $.

Let $\ell=2(T+T^3)$ where $T$ is a maximal mixing time. For $t \geq
\ell
$, a set of particle pairs $A$ and $(x,y) \in A$, let $\mathcal
{B}_{(x,y)}(t)$ denote the following event: There is no $ab$ meeting
for any $(a,b) \in A$ in the period $[\ell,t-1]$, and only $xy$ meet at
time $t$. Then
%
\begin{equation}\qquad
\Pr \bigl(\cB_{(x,y)}(t) \bigr)= {\biggl(1+O{\biggl(\frac{1}{n^{\Omega(1)}}+
\frac{1}{(k
\ln n)^{\Omega^+(1)}}+\frac{k^4t}{n^2}\biggr)}\biggr)}p\bigl(1-|A|p\bigr)^t,
\label{probBshorthand}
\end{equation}
where
%
\begin{equation}
p=\frac{1}{\theta_r n}{\biggl(1+O{\biggl(\frac{1}{n^{\Omega
(1)}}+\frac
{1}{(k \ln n)^{\Omega^+(1)}}
\biggr)}\biggr)}.\label{longformp}
\end{equation}
\end{lemma}

\begin{pf}
Let $S \subset V_H$ be the set of vertices in $H$ which correspond to
$x$ and $y$ being incident in $G$, but no other pair $(a,b) \in A$
being incident.
$|S|\leq n^{k-1}$ and $|S| \geq n^{k-1}-|A|^2n^{k-2} = n^{k-1}
{(1-O(\frac{k^4}{n}))}$. By similar calculations as in Lemma~\ref{pALem},
we get
\[
p_{v}=\frac{1}{\theta_r n}{\biggl(1+O{\biggl(\frac{1}{n^{\Omega
(1)}}+
\frac
{1}{(k \ln n)^{\Omega^+(1)}}\biggr)}\biggr)},
\]
where $v=\gamma(S)$.

We use Lemma~\ref{whichvx} with Lemma~\ref{pALem}. Referring to
\eqref
{probB}, $L\pi_{\gamma }=O(k^5(\ln n)^3/n)$ and
\[
\epsilon _v=O{\biggl(\frac{k^2}{n^{\Omega(1)}}+\frac{k^2}{(k \ln
n)^{\Omega^+(1)}}\biggr)}
=O{\biggl( \frac{1}{n^{\Omega(1)}}+\frac{1}{(k \ln n)^{\Omega^+(1)}}\biggr)}.
\]
Subsequently, \eqref{probB} gives us \eqref{probBshorthand}.
\end{pf}

Thus far, we have presented a set of lemmas that allow us to calculate
probabilities for meetings between subset of the $k$ particles. This
was done via the product graph framework, that allowed us to map the
state of the $k$ walks on $G$ to a single walk on $H$, then analyse
meetings between walks on $G$ in terms of the single walk on $H$
visiting specific vertices. From here on, we can largely forget about
$H$, and just use the lemmas we have established through it, in
particular, Lemma~\ref{ProbBxyLemma}, which will be the main tool used
to calculate probabilities of outcomes of the process.

In the next section, we will describe the interaction graph framework,
which allows us to map the unfolding of the process into a set of edge
weights that capture timing and outbreak information. We will then use
the tools from this section to calculate probabilities of particular
interaction graphs being realised.


\section{\texorpdfstring{Interaction graph.}{Interaction graph}}\label{Interaction graph}
The \emph{interaction graph} $\Upsi=(\mathcal{P}, E_{\Upsi})$ is a
weighted complete graph on the particle set $\mathcal{P}$, thus
$E_{\Upsi} = \{(x,y) \dvtx x,y \in\mathcal{P}, x \neq y\}$. For a particle~$x$, let $t(x)$ be the time at which $x$ is infected, or $\infty$ if it
never gets infected. For an \emph{interaction edge} $e=(x,y)\in
E_{\Upsi
}$, let $t(e)=\min\{t(x), t(y)\}$ [meaning $t(e) = \infty$ if neither
$x$ nor $y$ gets infected]. Then the weight $w_{\Upsi}(e)$ of the edge
is a random variable defined as
%
\begin{equation}\label{intGraphWeightDef}
w_{\Upsi}(e) = \cases{ %
\min \bigl
\{t-t(e) \dvtx t>t(e), \mbox{ $xy$ interaction at $t$} \bigr\}, &\quad
$\mbox{if $t(e)<
\infty$},$
\vspace*{2pt}\cr
\infty,& \quad $\mbox{otherwise}.$}\hspace*{-30pt}
\end{equation}

In the particular case that $\Row=1$, an interaction happens with every
meeting, in which case the edge weight represents the time elapsed
between the infection and the next meeting.

For the SI model, that is, $\xi=\infty$, every particle will get
infected almost surely, and so no edge of $\Upsi$ will have weight
$\infty$. This will not, in general, be the case for the SIR model.

We can think of a timer or clock associated with each interaction edge
$e=(x,y)$. The clock becomes active when either $x$ or $y$ becomes
infected, and stops when they next interact. Hence, when a particle $x$
gets infected, one clock stops being active, and generally, some number
become active simultaneously. The exception to the latter is when all
other particles became infected before $x$, in which case, the clocks
associated with $x$ had already become active previously.

Thus, the weighted complete graph $\Upsi$ can be represented by a
random ${k\choose2}$-dimensional vector $(w_{\Upsi}(e_i))_{i=1,\ldots,
{k\choose2}}$ where the edges are labelled in some arbitrary order. A
realisation of the graph $\Upsi$ is a specific set of values $\mathbf
{z}=(z_1, z_2, \ldots, z_{{k\choose2}})$ for each random variable
$w_{\Upsi}(e_i)$.


\section{\texorpdfstring{The SI model with $\rho=1$.}{The SI model with $rho=1$}}\label{SIModel}
This section deals with the special case $\xi=\infty$, $\Row=1$.

Two key ideas we use are that (1) the meeting times between pairs of
particles are almost independent, and that (2) the meeting time for a
pair of particles $(x,y)$ is roughly distributed as $\Ge(\frac
{1}{\theta
_rn})$, that is, as the geometric distribution with parameter $\frac
{1}{\theta_rn}$. We formalise (1) and (2) in Section~\ref{interapprox}.

Let $F$ be a weighted graph on $\mathcal{P}$. For a particle $x \in
\mathcal{P}$, denote by $d_F(x)$ the weighted distance (i.e., shortest
weighted path length) from the initial infective $x_0$ to $x$.
Furthermore, for an edge (i.e., particle pair) $e=(x,y)$, let
$d_F(e)=\min\{d_F(x), d_F(y)\}$. The interaction graph allows us to
relate weighted distance to infection time. Lemma~\ref{timeDistance}
below holds for $\Row\leq1$.

\begin{lemma}\label{timeDistance}
For a particle $x$, $t(x)=d_{\Upsi}(x)$, and for an edge $e \in
E_{\Upsi
}$, $t(e)= d_{\Upsi}(e)$.
\end{lemma}

\begin{pf}
For a particle $x$, let $P=(p_0=x_0,p_1, p_2, \ldots, p_l=x)$ be a
shortest path from $x_0$ to $x$ in $\Upsi$. For $0 \leq i \leq l-1$,
observe $t(p_{i+1}) \leq t(p_i)+w_{\Upsi}(p_i, p_{i+1})$, thus,
iterating, $t(x)\leq\sum_{i=0}^{l-1}w_{\Upsi}(p_i, p_{i+1})=d_\Upsi(x)$.

Now, for $x$ to get infected, there must be a chain of infections
$Q=(x_0=q_0, q_1, \ldots, q_{l'}=x)$ from one particle to another
starting from $x_0$.
For $0 \leq j \leq l'-1$, since $q_{j+1}$ is infected by $q_j$, we have
$t(q_{j+1})=t(q_j)+w_{\Upsi}(q_j,q_{j+1})$. Thus, iterating, we have
$t(x) = \sum_{j=0}^{l'-1}w_{\Upsi}(q_j,q_{j+1}) \geq d_{\Upsi}(x)$.
\end{pf}

In subsequent sections, this reduction is used to determine a
completion time for the process and to derive another graph that gives
the number of infected particles in the case $\xi<\infty$.

\subsection{\texorpdfstring{(Almost) building an interaction graph.}{(Almost) building an interaction graph}}\label{UpsilonPrimeSection}
One can build $\Upsi$ simply by observing the process unfolding and
starting and stopping clocks as it does so to determine edge weights.
To \emph{calculate} the probability of $\Upsi$ taking a particular
value [represented as a ${k\choose2}$-vector of edge weights], we will
use Lemma~\ref{ProbBxyLemma}. However, this lemma only allows us to
calculate the meeting times at values of $t \geq\ell=2(T+T^3)$ where
$T=O(k \ln n)$ is the mixing time. As such, we couple the process to a
slightly modified version of it, and construct the interaction graph
under the new process. We shall refer to this interaction graph as
$\Upsi'$.

Recall the role of clocks in determining edge weights of $\Upsi$; the
clock associated with $e=(x,y)$, $\mathit{clock}(e)$, is active precisely in the
period $[t(e), t'-1]$, where $t(e)=\min\{t(x), t(y)\}$, and $t'$ is the
first time after $t(e)$ that $x$ and $y$ interact. We describe the edge
$e$ as being \emph{active} while $\mathit{clock}(e)$ is active. Thus there is a
sequence of at most ${k\choose2}$ times $\t_0, \t_1, \t_2, \ldots,
\t
_j$, in which the set of active edges changes. We call these times
\emph
{epochs}. We let $\t_0=0$, since at this point, $x_0$ has active edges
with each of the other $k-1$ particles.

We parameterise the interaction probability so that $\Row(t)$ is the
probability of interaction at time $t$, universally for all particle
pairs. The new process makes the following modifications to the
original one:

\begin{longlist}[(ii)]
\item[(i)] {set $\Row(t)=0$ for $t \in[1, \ell]$;}
\item[(ii)] {if at time $\t$ there was an interaction between an
active pair $(x,y)$, set $\Row(t)=0$ for $t \in[\t+1, \t+\ell]$.}
\end{longlist}

These are the only differences; at all other times, we keep $\Row
(t)=1$. Furthermore, the definition of \emph{active edge} remains the
same for $\Upsi'$; \emph{an edge (particle pair) $e=(x,y)$ is active
precisely in the period $[t(e), t'-1]$, where $t(e)=\min\{t(x),t(y)\}$
and $t'$ is the first time after $t(e)$ that $x$ and $y$ interact.}

In the periods where $\Row(t)=0$, the set of active edges does not
change, and infections are not passed on, even if there were meetings
between infected and susceptible particles. We call these \emph{blind}
periods, because setting $\Row(t)=0$ is like ignoring interactions that
may have otherwise occurred during those times. We can think of $\ell$
as an extended mixing time.

Now we define a weighted complete graph $\Upsi'=(\mathcal{P},
E_{\Upsi
'})$: this is exactly $\Upsi$ but with edge weights determined in the
modified process. That is, for an edge $e=(x,y)$, the weight $w_{\Upsi
'}(e)$ of the edge is a random variable defined by the RHS of \eqref
{intGraphWeightDef}, but now $t(x)$ and $t(y)$ are the infection times
assuming blind periods. Observe $w_{\Upsi'}(e) > \ell$ because we still
start the clock associated with it at time $t(e)$, but due to the
blind periods, we will not stop it until the first time after
$t(e)+\ell
$ that they interact.

The modifications mean that in the sequence of epochs $\t_0, \t_1, \t
_2, \ldots, \t_j$ in which the set of active edges changes, we have
$\t
_{i+1}-\t_i>\ell$, and for each $i \in\{0,1,\ldots,j\}$ we set
$\Row
(t)=0$ for $t \in[\t_i+1,\t_i+\ell]$.

The initial stages of the construction of $\Upsi'$ for the SI model
with $\Row=1$ is described as follows: We set $\Row(t)=0$ for the first
$\ell$ steps and then set $\Row(t)=1$, waiting until $x_0$ interacts
with some particle. Suppose this happens with $x_1$ at time $\tau_1$.
We label the edge $(x_0, x_1) \in\Upsi'$ with $\tau_1$. We then
$\Row
(t)=0$ at steps $\tau_1+1, \ldots, (\tau_1+\ell)$, then resume
$\Row
(t)=1$, waiting for the next interaction to take place, either between
$x_0$ and one of the remaining $k-2$ particles, or $x_1$ and the
remaining $k-2$ particles. Suppose $x_2$ is the next amongst the
remaining particles to interact, and this happens at time $\t_2$. If it
happens with $x_0$, then the edge $(x_0, x_1) \in\Upsi'$ has weight
$\tau_2$. If it is with $x_1$, then the edge $(x_1,x_2) \in\Upsi'$ has
weight $\tau_2-\tau_1$. We then $\Row(t)=0$ for the following $\ell$
steps then resume $\Row(t)=1$ thereafter, waiting for the next relevant
interaction.

Thus, the partial construction of $\Upsi'$ at $\t_2$ in the former case
has $\Upsi'=\{w_{\Upsi'}(x_0, x_1)=\t_1, w_{\Upsi'}(x_0, x_2)=\t
_2\}$,
and in the latter case, it is $\Upsi'=\{w_{\Upsi'}(x_0,\break x_1)=\t_1,
w_{\Upsi'}(x_1, x_2)=\t_2-\t_1\}$. At time $\t_2$, in the former case,
$x_0$ has active clocks with $k-3$ other particles, and each of $x_1$
and $x_2$ with $k-2$ others (including each other). Suppose in this
scenario, that the next interaction after $\t_2+\ell$ occurs between
$x_1$ and $x_2$, at $\t_3$. Then $\Upsi'$ at this point will be
$\Upsi
'=\{w_{\Upsi'}(x_0, x_1)=\t_1, w_{\Upsi'}(x_0, x_2)=\t_2, w_{\Upsi
'}(x_1, x_2)=\t_3-\t_1\}$.

Continuing in this manner, we will eventually build the complete
edge-weighted graph $\Upsi'$. \emph{This will be exactly the same as
$\Upsi$ if, in the original process, there were no active pair
interactions when $\Row(t)=0$ in the modified process}. If, in fact,
there were active pair interactions at those time steps, then $\Upsi$
will be different to $\Upsi'$ described above. However, we will show
that, w.h.p., this will not be the case. That is, w.h.p., the graph
$\Upsi
'$, we construct in this modified process, is the same as $\Upsi$ that
would have been constructed had we kept $\Row(t)=1$ throughout.

\subsection{\texorpdfstring{Interaction graph approximation.}{Interaction graph approximation}}\label{interapprox}
We define a complete graph $\Lambda=(V_{\Lambda}, E_{\Lambda})$ on the
particle set $\mathcal{P}$ with i.i.d. random edge weights. The weight
$w_{\Lambda}(e)$ of each edge $e\in E_{\Lambda}$ has distribution
$\Ge
(q)$, where $q = \frac{\psi}{\theta_r n}$. In this section, $\psi=1$
since $\Row=1$.

$\Lambda$ can be specified by a ${k\choose2}$-vector; labelling the
edges in $\Lambda$ as $e_i$ with $1\leq i \leq{k\choose2}$, let $z_i$
denote the realised weight of $e_i$, and let $\Lambda=\mathbf{z}$ where
$\mathbf{z}=(z_1, z_2, \ldots, z_{{k\choose2}})$ denote this particular
realisation of $\Lambda$. Then
%
\begin{equation}
\Pr(\Lambda=\mathbf{z})=\Pr \Biggl(\bigwedge_{i=1}^{{k\choose
2}}w_{\Lambda}(e_i)=z_i
\Biggr) = \prod_{i=1}^{{k\choose
2}}q(1-q)^{z_i-1}\label{dualprob}.
\end{equation}

Our strategy will be as follows: We shall calculate the probability
$\Pr
(\Upsi'=\mathbf{z})$ of a particular realisation $\mathbf{z}$ of
$\Upsi
'$ by repeated application of Lemma~\ref{ProbBxyLemma}. We shall show
that this probability is well approximated by $\Pr(\Lambda=\mathbf
{z})$; that is, $\Lambda$ serves as an ``idealised'' version of $\Upsi
'$ (and, in turn, $\Upsi)$. Hence, w.h.p. results in $\Lambda$ can be
transferred to $\Upsi'$, and in turn transferred to $\Upsi$ if the
latter two are the same. Thus, to complete this proof strategy, we will
show that $\Upsi$ and $\Upsi'$ are the same w.h.p.

Assume the edges of $\Upsi'$ and $\Lambda$ are labelled in the same
order. We shall show that for a class of $\mathbf{z}$ defined below as
\emph{good}, $\Pr(\Lambda=\mathbf{z})$ is a close approximation for
$\Pr(\Upsi'=\mathbf{z})$.

\begin{definition}[({Good})]\label{goodDef}
Let $T$ be a maximal mixing time and let $\ell= 2(T+T^3)$. We say
$\mathbf{z}$ is \emph{good} when all of the components of $\mathbf{z}$
are finite and it satisfies the following:
\begin{longlist}
\item[(a)] if $\Upsi=\mathbf{z}$, then none of the
interactions that form the edges of $\Upsi$ occur within $\ell$ steps
of each other;
\item[(b)] $\sum_{i=1}^{{k\choose2}}z_i \leq k^2n \ln n$.
\end{longlist}
We shall also refer to a graph $F$ as good if $F=\mathbf{z}$ and
$\mathbf{z}$ is good.
\end{definition}

Part {(a)} of the above implies that when an infection occurs,
no other infection takes place at that time step, nor during the
following $\ell$. We require this condition because we wish to apply
Lemma~\ref{ProbBxyLemma}, which only gives probabilities for meetings
that occur after $\ell$ steps. Therefore, we can only use Lemma~\ref
{ProbBxyLemma} to calculate probabilities for interaction graphs in
which all relevant interactions are separated by $\ell$ steps.

Part {(b)} is a technical condition. We require it because if
weights $z_i$ are too big, then our approximations will cease to hold.

As will be seen in Lemma~\ref{UpsiesEqual}, the probability of $\Upsi$
having weights violating this condition goes to zero asymptotically.

We remind that in the context of $\Upsi'$, the notation $t(x)$ for a
particle $x$ refers to the time at which $x$ is infected when we set
$\Row(t)=0$ in blind periods. Similarly, for an edge (particle pair)
$e=(x,y)$, $t(e) =\min\{t(x),t(y)\}$. We may set $t(x)= \infty$ for a
particle that never gets infected. In the present case where the
infectious period $\xi=\infty$ (i.e., an SI model), every particle will
get infected almost surely, but when we address $\xi<\infty$, this may
not be the case, and the convention of setting ``infection times'' to
$\infty$ will be convenient.

We shall use graph notation with vectors $\mathbf{z}$: For an edge
$e_i$, $w_\mathbf{z}(e_i)=z_i$; for a particle $y$, $d_{\mathbf{z}}(y)$
is the weighted distance under {z} of $y$ from $x_0$; for an
edge $e=(x,y)$, $d_{\mathbf{z}}(e) = d_{\mathbf{z}}(x,y)=\min\{
d_{\mathbf{z}}(x),d_{\mathbf{z}}(y)\}$.

Lemma~\ref{distanceTimeEquiv} below holds for the SI model; when we
consider the case $\xi<\infty$ we will give a generalisation of it.

\begin{lemma}\label{distanceTimeEquiv}
For good $\mathbf{z}$, consider the following:
\begin{longlist}[(ii)]
\item[(i)] $\Upsi'=\mathbf{z}$;
\item[(ii)] \emph{for each particle pair (edge) $e=(x,y)$, there is no
$xy$ interaction in the period, $[d_{\mathbf{z}}(e)+\ell, d_{\mathbf
{z}}(e)+w_{\mathbf{z}}(e)-1]$, and there is an $xy$ interaction at time
$\tau(e)=d_{\mathbf{z}}(e)+w_{\mathbf{z}}(e)$}.
\end{longlist}
Then \textup{{(i)}} holds if and only if \textup{{(ii)}} holds.
\end{lemma}

\begin{pf}
$\Upsi' =\mathbf{z}$ can be restated as: For each $e=(x,y)$, the first
$xy$ interaction after time $t(e)+\ell$ is at time $t(e)+w_\mathbf
{z}(e)$. Hence, if we show that for each particle $x$, $t(x)=d_{\mathbf
{z}}(x)$, then we are done.

{(i)${}\Rightarrow{}$(ii)}
An equivalent to Lemma~\ref{timeDistance} holds for $\Upsi'$. The
implication follows.

{(ii)${}\Rightarrow{}$(i)}
Order particles by their distance from $x_0$ in $\mathbf{z}$:
$x_0=x_{(0)}, x_{(1)}, \ldots,\break x_{(k-1)}$ where $i<j \Rightarrow
d_{\mathbf{z}}(x_{(i)}) \leq d_{\mathbf{z}}(x_{(j)})$. We shall prove
that $t(x_{(i)})=d_{\mathbf{z}}(x_{(i)})$.

Clearly this proposition holds for $x_{(0)}$. Suppose for all $i \leq
N-1$, $t(x_{(i)})=d_{\mathbf{z}}(x_{(i)})$. If $N=k$, we are done.
Otherwise $N<k$. Let $x_{(M)}$ be a neighbour of $x_{(N)}$ on a
shortest path from $x_0$ to $x_{(N)}$. Since $d_{\mathbf{z}}(x_{(M)}) <
d_{\mathbf{z}}(x_{(N)})$, $M < N$, so by the induction hypothesis,
$t(x_{(M)})=d_{\mathbf{z}}(x_{(M)})=d_{\mathbf{z}}(e)$ where
$e=(x_{(M)}, x_{(N)})$. By {(ii)} this implies $t(x_{(N)}) \leq
d_{\mathbf{z}}(e)+w_\mathbf{z}(e)=d_{\mathbf{z}}(x_{(N)})$.

Now consider the chain of infections starting at $x_0$ that led to
$x_{(N)}$ being infected. Let $x_{(j)}$ be the first in the chain where
$j>N-1$ and suppose it got infected by $x_{(i)}$, $i \leq N-1$. Then by
{(ii)}, $t(x_{(j)}) = d_{\mathbf{z}}(x_{(i)})+w_{\mathbf
{z}}(x_{(i)}, x_{(j)})\geq d_{\mathbf{z}}(x_{(j)}) \geq d_{\mathbf
{z}}(x_{(N)})$, implying $t(x_{(N)}) \geq d_{\mathbf{z}}(x_{(N)})$.

Hence $t(x_{(N)})=d_{\mathbf{z}}(x_{(N)})$ and the lemma follows.
\end{pf}

Lemma~\ref{ineractionapprox} below holds for the SI model; when we
consider the case $\xi<\infty$, we will give a generalisation of it.

\begin{lemma}\label{ineractionapprox}
Assume the conditions of Section~\ref{assumptionssection}.

For good $\mathbf{z}$,
%
\begin{equation}
\Pr \bigl(\Upsi'=\mathbf{z} \bigr) = \bigl(1+o(1) \bigr)\Pr(
\Lambda= \mathbf{z}). \label
{equivprobs}
\end{equation}
\end{lemma}

\begin{pf}
By Lemma~\ref{distanceTimeEquiv}, $\Upsi'=\mathbf{z}$ defines, for each
edge $e=(x,y)$, a time $\tau(e)=d_{\mathbf{z}}(e)+w_{\mathbf{z}}(e)$.
Letting $\t_0=t(x_0)=0$, since $\mathbf{z}$ is good, we get a sequence
of epochs $\tau_0<\tau_1< \tau_2< \cdots<\tau_{{k\choose2}}$ that are
at least $\ell$ apart.

An edge $e$ is active precisely in the period $[d_{\mathbf{z}}(e),
d_{\mathbf{z}}(e)+w_{\mathbf{z}}(e)-1]$ (and at no other time). Let
$A_i$ denote the set of active edges in the period $[\t_i, \t
_{i+1}-1]$. This defines a sequence $(A_i)=(A_0, A_1, \ldots,
A_{{k\choose2}-1})$ of active edge sets associated with epochs $\tau
_0, \tau_1, \ldots, \tau_{{k\choose2}-1}$, respectively, and which
remain constant between epochs.
The active set changes at each epoch, when an active edge is removed
from the set, and possibly new ones are added. In general, an edge $e$
will be a member of a number of active edge sets $A_i$.

Let $\sigma(\mathbf{z})= ((e_{(1)},\tau_1), (e_{(2)}, \tau_2),
\ldots, (e_{ \bigl({k\choose2} \bigr)}, \tau_{{k\choose2}})
)$ be
defined by the above, where $e_{(i)}$ is the particle pair that
interact at time $\t_i$. We shall use $\sigma(\mathbf{z})$ and $(A_i)$
to calculate $\Pr(\Upsi'=\mathbf{z})$. In particular, $\Upsi
'=\mathbf
{z}$ if and only if both of the following hold:
\begin{longlist}[(1)]
\item[(1)] if $(x,y) \in A_i$, there is no $xy$ interaction in the
period $[\t_i+\ell, \t_{i+1}-1]$;
\item[(2)] the particle pair $e_{(i)}$ interact at time $\t_i$, and
this is the only pair in $A_i$ that does.
\end{longlist}

The probability of the above will be determined by repeated application
of Lemma~\ref{ProbBxyLemma} and the strong Markov property. Consider
the process up until epoch~$\t_1$; {(1)} dictates no interaction
between any $(x_0,y) \in A_0=\{(x_0,y)\dvtx y \in\mathcal{P}, y \neq x_0\}
$ in the period $[\ell, \t_1]$, and {(2)} dictates only an
$e_{(1)} = (x_0, x_1)$ interaction at $\tau_1$. Applying Lemma~\ref
{ProbBxyLemma}:
%
\begin{equation}\label{fdgdgdg8}
\Pr \bigl(\cB_{e_{(1)}}(\t_1) \bigr)= {\biggl(1+O{\biggl(
\frac{1}{n^{\Omega
(1)}}+\frac
{1}{(k \ln n)^{\Omega^+(1)}}+\frac{k^4\t_1}{n^2}\biggr)}
\biggr)}p\bigl(1-|A_0|p\bigr)^{\t_1}.\hspace*{-25pt}
\end{equation}

Now if we consider the next period $[\tau_1, \tau_2-1]$, the set of
active edges are $A_1 = \{(x_0, y) \dvtx y \in\mathcal{P}, y \neq x_0,
x_1\} \cup\{(x_1, y) \dvtx y \in\mathcal{P}, y \neq x_0, x_1\}$.
{(1)} dictates no interaction between any pair in $A_1$, and
{(2)} dictates only an $e_{(2)}$ interaction at $\t_2$. Thus we can
apply Lemma~\ref{ProbBxyLemma} for this period and active edge set to get
\begin{eqnarray*}
&&\Pr \bigl(\cB_{e_{(2)}}(\t_2) \bigr)\\
&&\qquad= {\biggl(1+O{\biggl(
\frac{1}{n^{\Omega
(1)}}+\frac
{1}{(k \ln n)^{\Omega^+(1)}}+\frac{k^4(\t_2-\t
_1)}{n^2}\biggr)}
\biggr)}p\bigl(1-|A_1|p\bigr)^{\t
_2-\t_1}.
\end{eqnarray*}

Because $\mathbf{z}$ is good, $(\t_{j+1}-\t_j)$ is such that $k^4(\t
_{j+1}-\t_j)/n^2=n^{-\Omega(1)}$, which can be absorbed into the
correcting factor of $p$, which has form \eqref{longformp}. Therefore,
we can write $\Pr(\cB_{e_{(1)}}(\t_1))= p(1-|A_0|p)^{\t_1}$ and
$\Pr(\cB
_{e_{(2)}}(\t_2))=p(1-|A_1|p)^{\t_2-\t_1}$.

We can continue in similar fashion for each epoch, and by the strong
Markov property, we can multiply these probabilities to get
%
\begin{equation}
\Pr \bigl(\Upsi'=\mathbf{z} \bigr) = \prod
_{j=0}^{{k\choose
2}-1}p\bigl(1-|A_j|p\bigr)^{\t
_{j+1}-\t_j},
\label{faoifs}
\end{equation}
where
\[
p=\frac{1}{\theta_r n}{\biggl(1+O{\biggl(\frac{1}{n^{\Omega
(1)}}+\frac
{1}{(k \ln n)^{\Omega^+(1)}}
\biggr)}\biggr)}.
\]

Letting $\Delta_j=\tau_{j+1}-\tau_j$, \eqref{faoifs} can be written
\begin{eqnarray*}
\Pr \bigl(\Upsi'=\mathbf{z} \bigr) &=& p^{{k\choose2}}\prod
_{j=0}^{{k\choose
2}-1}\exp \bigl\{-
\bigl(1+O(|A_j|p) \bigr)|A_j|p\Delta_j
\bigr\}
\\
&=&p^{{k\choose2}}\exp \Biggl\{- \bigl(1+O \bigl(k^2p \bigr)
\bigr)p \sum_{j=0}^{{k\choose
2}-1}|A_j|
\Delta_j \Biggr\}.
\end{eqnarray*}
Furthermore,
\[
\sum_{j=0}^{{k\choose2}-1}|A_j|
\Delta_j=\sum_{j=0}^{{k\choose
2}-1}\sum
_{e}\mathbf{1}_{\{e \in A_j\}}\Delta_j
= \sum_{e}\sum_{j=0}^{{k\choose
2}-1}
\mathbf{1}_{\{e \in A_j\}}\Delta_j
\]
and
\[
\sum_{j=0}^{{k\choose2}-1}\mathbf{1}_{\{e \in A_j\}}
\Delta_j = \sum_{j=0}^{{k\choose2}-1}
\mathbf{1}_{\{e \in A_j\}}(\tau_{j+1}-\tau_j)=z_e
\]
because this sums over all intervals $[\tau_j, \tau_{j+1}-1]$ in which
edge $e$ is active, that sum being $z_e$. Hence, $\sum_{j=0}^{{k\choose
2}-1}|A_j|\Delta_j=\sum_{e}z_e = \sum_{i=1}^{{k\choose2}}z_i$.

Thus
\[
\Pr \bigl(\Upsi'=\mathbf{z} \bigr)=p^{{k\choose2}}\exp \Biggl\{ -
\bigl(1+O \bigl(k^2p \bigr) \bigr)p\sum_{i=1}^{{k\choose2}}z_i
\Biggr\}.
\]
From \eqref{dualprob}, we have
\[
\Pr(\Lambda=\mathbf{z})= \prod_{i=1}^{{k\choose2}}q(1-q)^{z_i-1}
= q^{{k\choose2}}(1-q)^{-{k\choose2}}\exp \Biggl\{- \bigl(1+O(q) \bigr)q\sum
_{i=1}^{{k\choose2}}z_i \Biggr\}.
\]
Therefore,
\begin{eqnarray}
\qquad&&\frac{\Pr(\Upsi'=\mathbf{z})}{\Pr(\Lambda=\mathbf{z})}\nonumber\\
&&\qquad= \biggl(\frac {p(1-q)} {q}\biggr)^{{k\choose2}}\exp \Biggl\{
{\bigl( \bigl(1+O(q) \bigr)q- \bigl(1+O \bigl(k^2p \bigr) \bigr)p
\bigr)}\sum_{i=1}^{{k\choose2}}z_i
\Biggr\},\nonumber
\\
\label{fsg4d56}
&&\biggl(\frac{p(1-q)} {q}^{{k\choose2}}\biggr)={\biggl(1+O{\biggl(
\frac
{1}{n^{\Omega
(1)}}+ \frac{1}{(k \ln n)^{\Omega^+(1)}}\biggr)}\biggr)}^{{k\choose2}}
\nonumber
\\[-8pt]
\\[-8pt]
\nonumber
&&\hspace*{48pt}\qquad =1+O{\biggl(\frac{1}{n^{\Omega(1)}}+\frac{1}{(k \ln n)^{\Omega
^+(1)}}\biggr)},
\end{eqnarray}
where we have used the fact that $k\leq n^{\epsilon}$ for $\epsilon$
small enough.

Since
\[
\bigl(1+O \bigl(k^2p \bigr) \bigr)p- \bigl(1+O(q) \bigr)q=O{\biggl(
\frac{1}{n^{\Omega(1)}}+\frac{1}{(k
\ln
n)^{\Omega^+(1)}}\biggr)}\frac{1}{n},
\]
then by part {(b)} of the definition of \emph{good} {z},
%
\begin{equation}
{\biggl(\frac{1}{n^{\Omega(1)}}+\frac{1}{(k \ln n)^{\Omega
^+(1)}}\biggr)}\frac
{1}{n}\sum
_{i=1}^{{k\choose2}}z_i <
\frac{k^2\ln n}{n^{\Omega
(1)}}+ \frac
{k^2\ln n}{(k \ln n)^{\Omega^+(1)}}. \label{l4l4j}
\end{equation}

Since the $\Omega^+(1)$ term in \eqref{fsg4d56} and \eqref{l4l4j} is an
arbitrarily large constant, they are $1+o(1)$ and $o(1)$, respectively,
when $k \leq n^\epsilon $ for sufficiently small $\epsilon $. Thus,
we conclude $\Pr
(\Upsi'=\mathbf{z})/\Pr(\Lambda=\mathbf{z})=1+o(1)$.
\end{pf}

To prove Lemma~\ref{UpsiesEqual} below, we require the following minor
adaptation of Lemma~20 in \cite{CFR-Mult}:

\begin{lemma}[(\cite{CFR-Mult})]\label{dontmeet}
Let $G$ be typical, and let $k \leq n^\epsilon$ for sufficiently small~$\epsilon$.

Suppose that particles start with minimum separation at least
$d=\break\a(\ln\ln n + \ln k)$
where $\alpha$ is a constant. Let $\tau=O(T^3)$ where $T=O(k \ln n)$. Then
%
\begin{equation}\label{84g64sdg}
\Pr(\mbox{a given pair of particles $x,y$ meet during $\tau $})=O \bigl(\tau
^2/(r-1)^{d/6} \bigr).\hspace*{-20pt}
\end{equation}
\end{lemma}

It should be noted that the $\alpha$ in Lemma~\ref{dontmeet} can be
made as large as required by making $\epsilon$ sufficiently small.

Lemma~\ref{UpsiesEqual} below holds for general $\xi$ and $\Row\leq1$.

\begin{lemma}\label{UpsiesEqual}
Let $G$ be typical, and let $k \leq n^\epsilon$ for sufficiently small
$\epsilon$.

With high probability:
\begin{longlist}[(a)]
\item[(a)] none of the interactions that form the finitely-weighted
edges of $\Upsi$ occur within $\ell$ steps of each other;
\item[(b)] the sum of the finite edge weights of $\Upsi$ is at most
$k^2n \ln n$.
\end{longlist}
\end{lemma}

\begin{pf}
It will be convenient to prove part {(b)} of the lemma first.

{(b)} Consider an active pair of particles $(x,y)$. The expected
number of meetings they have before interacting is $1/\Row$. The
expected time between meetings is $O(n)$. Since walks and interactions
are independent, we can multiply these quantities together to give a
bound on the expected time till an interaction. Since $\Row$ is assumed
to be constant, this is $O(n)$.

The result now follows by linearity of expectation and Markov's inequality.

{(a)} Let $A$ be a set of active particle pairs. We let them mix
for $\ell$ steps, and suppose the first interaction occurs at some
(random) time $t \geq\ell$. Let $S \subset V_H$ correspond to a
meeting of at least one pair in $A$. Let $S' \subset S$ correspond to
more than one pair being within distance $d$ of each other. We wish to
calculate $\Pr(\mathcal{W}_{\ul u}^H(t) \in S' \mid\mathcal{W}_{\ul
u}^H(t) \in S \wedge\mathcal{W}_{\ul u}^H(\t) \notin S \mbox{ for }
\t
\in[\ell, t-1])$.

Let $\gamma =\gamma (S)$ and $\gamma '=\gamma (S')$. By Lemma~\ref
{pALem}, $\p_\gamma
=(1+o(1)|A|/n$. Furthermore, since $|S'| \leq{|A|\choose
2}r^{2d}n^{k-2}$, $p_{\gamma '}=O(|A|^2r^{2d}/n^2)$.

Applying Lemma~\ref{whichvx},\footnote{Note, although the lemma is
stated for a visit to a vertex $v$, $v$ can also be a contraction of a
subset of $S$.}
%
\begin{eqnarray}\label{dsggh564g}
&&\Pr \bigl(\mathcal{W}_{\ul u}^H(t) \in S'
\wedge\mathcal{W}_{\ul u}^H(\t) \notin S \mbox{ for } \t\in[
\ell, t-1] \bigr)\nonumber
\\
&&\qquad= \biggl(1+O(L\pi_{\gamma })+O(\epsilon _{\gamma '})+O \biggl(
\frac
{k^4t}{n^2} \biggr) \biggr)p_{\gamma '} (1-p_{\gamma })^t
\\
&&\qquad= \bigl(1+o(1) \bigr)p_{\gamma '} (1-p_{\gamma })^t.
\nonumber
\end{eqnarray}
The correcting factor in the last line holds because $t \leq k^ 2 n \ln
n$ by part {(b)}, and it is straightforward to show $\epsilon
_{\gamma '}=o(1)$.
%
\begin{eqnarray}\label{8g8x4h}
&&\Pr \bigl(\mathcal{W}_{\ul u}^H(t) \in S \wedge
\mathcal{W}_{\ul u}^H(\t) \notin S \mbox{ for } \t\in[\ell,
t-1] \bigr)
\nonumber
\\[-8pt]
\\[-8pt]
\nonumber
&&\qquad=\bigl(1+o(1)\bigr)p_\gamma (1-p_\gamma )^t,
\end{eqnarray}
so dividing \eqref{dsggh564g} by \eqref{8g8x4h} gives
\begin{eqnarray*}
\operatorname{Pr} \bigl(\mathcal{W}_{\ul u}^H(t) \in S' \mid
\mathcal{W}_{\ul u}^H(t) \in S \wedge\mathcal{W}_{\ul u}^H(
\t) \notin S \mbox{ for } \t\in[\ell, t-1] \bigr) &=& O \biggl(\frac{|A|r^{2d}} {
n}
\biggr)
\\
&=&O\biggl(\frac{1} {n^{\Omega(1)}}\biggr).
\end{eqnarray*}

There are at most ${k\choose2}$ interactions determining the edge
weights of $\Upsi$; taking the union bound over all of them, this is
$o(1)$ if $\epsilon$ is small enough.

The RHS of \eqref{84g64sdg} is $O(1/(k \ln n)^{\Omega^+(1)})$.
Consequently, we can apply Lem\-ma~\ref{dontmeet} across all [at most
${k\choose2}$] particles pairs across all [at most ${k\choose2}$]
interactions. We thus conclude that w.h.p., at most one pair interact at
any time~and there are no interactions in any of the following length
$\ell$ blind periods.
\end{pf}

\begin{corollary}\label{UpsiesAREequal}
With high probability:
\begin{longlist}[(iii)]
\item[(i)] for general $\xi$, $\Upsi=\Upsi'$;
\item[(ii)] when $\xi=\infty$, $\Upsi$ and $\Upsi'$ are good;
\item[(iii)] when $\Row=1$, $\Lambda$ is good.
\end{longlist}
\end{corollary}

Part {(iii)} of Corollary~\ref{UpsiesAREequal} is a corollary of
Lemmas \ref{ineractionapprox} and~\ref{UpsiesEqual} together:
observe $1-o(1)=\sum_{\mathbf{z}\ \mathrm{good}}\Pr(\Upsi'=\mathbf{z})
=(1+o(1))\sum_{\mathbf{z}\ \mathrm{good}}\Pr(\Lambda=\mathbf{z})$.

In conjunction with Lemma~\ref{ineractionapprox}, Corollary~\ref
{UpsiesAREequal} will allow us to say that events in $\Lambda$ and
$\Upsi$ have roughly the same probability.

Events in $\Upsi$, $\Upsi'$, and $\Lambda$, are subsets of $\Omega_k$,
the set of all possible weightings on the ${k\choose2}$ edges. Thus
$\Omega_k$ is the set of all ${k\choose2}$-vectors with nonnegative
entries. An event $\mathcal{E} \subseteq\Omega_k$ in $\Upsi$ occurs if
and only if $\Upsi=\mathbf{z}$ where $\mathbf{z} \in\mathcal{E}$,
indicated by $\mathbf{1}_{\mathcal{E}}(\mathbf{z})$ and similarly with
$\Upsi'$ and $\Lambda$.

Recall $\Lambda$ depends on $n$. Lemma~\ref{whpEqual} below holds for
the SI model; when we consider the case $\xi<1$, we will give a
generalisation of it.

\begin{lemma}\label{whpEqual}
Assume the conditions in Section~\ref{assumptionssection}.

Suppose there is a sequence of events $(\mathcal{E}_n)_{n \geq1}$ and
a constant $p_c$ such that in $\Lambda$, $\Pr(\mathcal
{E}_n)\rightarrow
p_c$ as $n \rightarrow\infty$. Then in $\Upsi$, $\Pr(\mathcal
{E}_n)\rightarrow p_c$ as $n \rightarrow\infty$.
\end{lemma}

\begin{pf}
Let $\mathbf{1}_{\mathcal{E}_n}(\mathbf{z})$ be the indicator that
$\mathbf{z} \in\mathcal{E}_n$. By Lemma~\ref{ineractionapprox} and
{(ii)} and {(iii)} of Corollary~\ref{UpsiesAREequal},
\begin{eqnarray*}
\Pr(\mathcal{E}_n) &=& \sum_{\mathbf{z}}\Pr(
\Lambda=\mathbf {z})\mathbf {1}_{\mathcal{E}_n}(\mathbf{z})
\\
&=& o(1)+\sum_{\mathbf{z}\ \mathrm{good}}\Pr(\Lambda=\mathbf {z})\mathbf
{1}_{\mathcal{E}_n}(\mathbf{z})
\\
&=& o(1)+ \bigl(1-o(1) \bigr)\sum_{\mathbf{z}\ \mathrm{good}}\Pr \bigl(
\Upsi'=\mathbf {z} \bigr)\mathbf{1}_{\mathcal{E}_n}(\mathbf{z})
\\
&=& o(1)+ \bigl(1-o(1) \bigr)\sum_{\mathbf{z}}\Pr \bigl(
\Upsi'=\mathbf{z} \bigr)\mathbf {1}_{\mathcal{E}_n}(\mathbf{z}).
\end{eqnarray*}
Hence in $\Upsi'$, $\Pr(\mathcal{E}_n) \rightarrow p_c$ as $n
\rightarrow\infty$. Consequently, by part {(i)} of Corollary~\ref{UpsiesAREequal}, $\Pr(\mathcal{E}_n) \rightarrow p_c$ as $n
\rightarrow\infty$ in $\Upsi$.
\end{pf}

\subsection{\texorpdfstring{Completion time.}{Completion time}}\label{comptimesection}
In \cite{CFR-Mult} an expectation of $\frac{2\theta_r n}{k}\ln k$ was
determined for the completion time of a broadcasting model on $k$
particles that is equivalent to the SI model with $\Row=1$. In Section~\ref{TheGeneralCase} the role of $\psi$ will be made clear.
Subsequently, the generalisation of the expectation to $\frac{2\theta_r
n}{\psi k}\ln k$ for $\Row\leq1$ will be seen to be straightforward.
In this section, we shall get a convergence in probability to the same
value, for the case $\Row=\psi=1$. Subsequent to the treatment in
Section~\ref{TheGeneralCase}, it will be clear how the result extends
to the general case.\looseness=-1

We make use of a theorem from \cite{Janson}: Assign each edge $(i,j)$
of a complete graph on $k$ vertices a random weight $Y_{ij}$. The
weights are assumed to be independent and identically distributed,
nonnegative and satisfying $\Pr(Y_{ij} \leq t) = t +o(t)$ as $t
\rightarrow0$. Let $X_{ij}$ be the minimal total weight of a path
between a given pair of vertices $i,j$.

\begin{theorem}[(\cite{Janson})]\label{Jansonsthm}
Under the assumptions above, for any fixed $i$, as $k \rightarrow
\infty$,
%
\begin{equation}
\frac{\max_{j}X_{ij}}{\ln k/k} \mathop{\xrightarrow}^{p} 2. \label{Jansoneq}
\end{equation}
\end{theorem}

We apply this below.

\begin{pf*}{Proof of Theorem~\ref{completiontimetheorem} for $\Row=1$}
If a random variable $Y \sim\mbox{Exp}(\lambda)$, that is, has
exponential distribution with parameter $\lambda$, then $\lambda Y
\sim
\mbox{Exp}(1)$, which is valid for Theorem~\ref{Jansonsthm}. Hence
consider a complete graph $F$ on the particle set with i.i.d. edge
weights $Y_{ij}\sim\mbox{Exp}(\lambda)$, and let $X_{ij}$ be the
minimal total weight of a path between a given pair of vertices
(particles) $i,j$ in $F$. By the above, as $k \rightarrow\infty$,
%
\begin{equation}
\frac{\max_{j}\lambda X_{x_0j}}{\ln k/k} \mathop{\xrightarrow}^{p} 2.
\end{equation}

Let $\lambda=1/(\theta_r n)$, so $w_\Lambda(i,j) \sim\Ge(\lambda
)$ for
an edge $e=(i,j)$ in $\Lambda$. For each particle pair $e=(i,j)$, let
$U_{ij}$ be i.i.d. random variables uniform on $[0,1]$. We use $U_{ij}$
to determine both $Y_{ij}$ in $F$ and $w_\Lambda(i,j)$ in $\Lambda$.
Specifically, let $Y_{ij}=\frac{-\ln(U_{ij})}{\lambda}$, and let
$w_{\Lambda}(i,j)=  \lceil\frac{\ln(U_{ij})}{\ln(1-\lambda)}
\rceil=  \lceil\frac{-\ln(U_{ij})}{(1+O(1/n))\lambda}
\rceil$. Hence $|Y_{ij}-w_{\Lambda}(i,j)| \leq1$ for all pairs
$i,j$ when $n$ is large enough, in which case, $|X_{x_0j}-d_\Lambda
(j)|<k^2$. Thus
\[
\frac{\max_{j}\lambda(X_{x_0j}-k^2)}{\ln k/k} \leq\frac{\max_{j}\lambda d_\Lambda(j)}{\ln k/k} \leq\frac{\max_{j}\lambda
(X_{x_0j}+k^2)}{\ln k/k},
\]
but $\lambda k^3/\ln k=O(k^3/n)$. Hence, since $k \leq n^\epsilon$,
when $\epsilon$ is small enough,
\[
\frac{\max_{j} d_\Lambda(j)}{n\ln k/k} \mathop{\xrightarrow}^{p} 2\theta_r.
\]

By Lemmas \ref{whpEqual} and \ref{timeDistance}, the lemma follows.
\end{pf*}
%

\section{\texorpdfstring{The SI(R) model with $\rho=1$.}{The SI(R) model with $rho=1$}}\label{SIRRho1}
Now we allow $\xi< \infty$, thus generalising the previous section. We
create $\Upsi$ and $\Upsi'$ exactly as before, and we will get
particular realisations of each. As before, for a particular
realisation $\Upsi=\mathbf{z}$ (or $\Upsi'=\mathbf{z}$), $\mathbf{z}$
is a ${k\choose2}$-vector; however, recalling the definitions from
Section~\ref{Interaction graph}, we see that now some of the entries
may be $\infty$, whereas before in the SI case, they were all finite
almost surely.
\ignore{
Note, for a particular realisation $\Upsi=\mathbf{z}$ (or $\Upsi
'=\mathbf{z}$), $\mathbf{z}$ is a ${k\choose2}$-vector as before,
except that some of the entries may be $\infty$.
}

We define the function $f_\xi(F)$ that takes as input a weighted graph
$F$ and returns the same graph but with all edges with weights
exceeding $\xi$ deleted.
We may also write $f_{\xi}(\mathbf{z})$, interpreting the argument as
some graph weighted by $\mathbf{z}$.

For a graph $F=(\mathcal{P}, E_F)$, denote by $\mathcal{C}_F$ the
connected component in $F$ that $x_0$ belongs to. For example, take a
complete graph on $\mathcal{P}$ weighted with $\mathbf{z}$, and apply
$f_{\xi}(\mathbf{z})$ to delete edges with weight exceeding $\xi$. The
connected component that $x_0$ belongs to in this graph is denoted by
$\mathcal{C}_{f_{\xi}(\mathbf{z})}$.

The following lemma tells us which particles ever get infected.

\begin{lemma}\label{PSiCompntInfect}
Let $F \in\{\Upsi, \Upsi'\}$. A particle $y$ becomes infected if and
only if $y \in\mathcal{C}_{f_\xi(F)}$.
\end{lemma}

\begin{pf}
If $y \in\mathcal{C}_{f_\xi(F)}$, then there is a path from $x_0$ to
$y$ with each edge having weight at most $\xi$. Each such edge can only
exist because one of the ends was infected, and they subsequently
interacted within $\xi$ steps after the infection, meaning that the
other will be infected after the interaction. Hence, $y$ got infected.

Conversely, suppose $y$ got infected. Then there was a chain of
infections from $x_0$ to $y$. This would define a path with each edge
weight at most $\xi$. Hence it would be in $\mathcal{C}_{f_\xi(F)}$.
\end{pf}

By Lemma~\ref{PSiCompntInfect}, we can determine the size of the
outbreak via $\mathcal{C}_{f_\xi(\Upsi)}$ (of course, for the SI case,
we did not need to consider $\mathcal{C}_{f_\xi(\Upsi)}$ since all edge
weights would be finite and so none would get deleted). As before, we
use the blind periods in order to apply the tools of Section~\ref{approachSection}. Thus we work with $\Upsi'$.

For a weighted graph $F=(\mathcal{P}, E_F)$ and a partition $A,B$ of
$\mathcal{P}$, denote by $E_F(A)$ and $E_F(B)$ the (weighted) edges
induced by $A$ and $B$, respectively. Denote by $E(A:B)$ the (weighted)
edges with one end in $A$ and the other in $B$. The weights are
explicitly part of the notation. Thus, for example, if we write
$E_{\Upsi'}(A)=E_{\mathbf{z}}(A)$ for some subset $A$ of the particles,
we equate both the set of edges induced, as well as their weights.

In $\Upsi'$, let $A$ and $B$ be the infected and noninfected
particles, respectively. Observe that all weights in $E_{\Upsi'}(B)$
are infinite, all in $E_{\Upsi'}(A:B)$ are finite but larger than $\xi
$, and all weights in $E_{\Upsi'}(A)$ are finite, with some being at
most $\xi$, and possibly some being greater than $\xi$.

The key observation is that \emph{given edge weights $E_{\Upsi'}(A)$
and $E_{\Upsi'}(A:B)$, we know edge weights $E_{\Upsi'}(B)$}. This is
because we have timing information on the edges that allow us to
reconstruct the unfolding of the process. Thus, determining the
probability in $\Upsi'$ of edge weights $E_{\Upsi'}(A)$ and $E_{\Upsi
'}(A:B)$ is all that it required to determine the rest of the system,
in particular, to show that those in $A$ get infected and those in $B$
do not.

Recall the definition of $\Lambda$ from Section~\ref{interapprox}. We
reason as follows: Draw a $\Lambda$ and delete all edges $e$ with
$w_\Lambda(e)>\xi$. The resulting component $\mathcal{C}_{f_\xi
(\Lambda
)}$ closely approximates $\mathcal{C}_{f_\xi(\Upsi)}$, informally
justified as follows: Call a path an \emph{$\xi$-path} if all edges in
it have weight at most $\xi$. Partition $\mathcal{P}$ into $A,B$ where
$y \in A$ iff there is a $\xi$-path from $x_0$ to $y$ in $\Lambda$.
Observe in $E_\Lambda(A:B)$ all weights are greater than $\xi$, and in
$E_\Lambda(A)$ there will be weights at most $\xi$, and possibly some
exceeding $\xi$. All weights in these two sets will be finite.
$E_{\Upsi
'}(A)$ and $E_{\Upsi'}(A:B)$ will have roughly the same probability
distribution on edge wights as $E_{\Lambda}(A)$ and $E_{\Lambda}(A:B)$,
and since $\Upsi=\Upsi'$ w.h.p., the same will be true of $\Upsi$.
Therefore, $\mathcal{C}_{f_\xi(\Upsi)}$ has roughly the same
probability distribution as $\mathcal{C}_{f_\xi(\Lambda)}$.
Consequently, probabilistic statements about the size of $\mathcal
{C}_{f_\xi(\Lambda)}$ can be transferred to $\mathcal{C}_{f_\xi
(\Upsi)}$.

We proceed to formalise the above. We require the following
generalisation of Lemma~\ref{distanceTimeEquiv} (proof in the
\hyperref[app]{Appendix}). We denote by $V(\mathcal{C}_{f_\xi(\mathbf{z})})$ the set of
vertices in $\mathcal{C}_{f_\xi(\mathbf{z})}$, that is, $V(\mathcal
{C}_{f_\xi(\mathbf{z})})=\{x \in\mathcal{P} \dvtx \mbox{ under
$\mathbf
{z}$ there is a $\xi$-path from}\break x_0\mbox{ to }x\}$.

\begin{lemma}\label{TimeDistPsi}
Suppose $\mathbf{z}$ is good, and let $A=V(\mathcal{C}_{f_\xi
(\mathbf
{z})})$, $B=\mathcal{P}\setminus A$. Consider the following:
\begin{longlist}[(ii)]
\item[(i)] $\Upsi'$ has $E_{\Upsi'}(A)=E_{\mathbf{z}}(A)$ and
$E_{\Upsi
'}(A:B)=E_{\mathbf{z}}(A:B)$;
\item[(ii)] for each particle pair (edge) $e=(x,y)$ such that
$d_{f_{\xi
}(\mathbf{z})}(e)<\infty$, there is no $xy$ interaction in the period,
$[d_{f_{\xi}(\mathbf{z})}(e)+\ell, d_{f_{\xi}(\mathbf
{z})}(e)+w_{\mathbf
{z}}(e)-1]$, and there is an $xy$ interaction at time $\tau
(e)=d_{f_{\xi
}(\mathbf{z})}(e)+w_{\mathbf{z}}(e)$.
\end{longlist}
Then {\textup{(i)}} holds if and only if \textup{{(ii)}} holds.
\end{lemma}

Thus, we have established that {(ii)} implies the set
$A=\mathcal
{C}_{f_\xi(\mathbf{z})}$ gets infected and that for each edge $e$ such
that $d_{f_\xi(\mathbf{z})}(e)<\infty$, that is, for each edge $e
\in
E_{\mathbf{z}}(A) \cup E_{\mathbf{z}}(A:B)$, $w_{\Upsi
'}(e)=w_{\mathbf
{z}}(e)$. It remains to establish that no particle in $B$ is infected.
If $d_{f_\xi(\mathbf{z})}(x)=\infty$, then $x$ could not have been
infected since if it was, then there must have been a chain of
infections from $x_0$ to $x$. In this chain, there must have been some
$x_b \notin\mathcal{C}_{f_\xi(\mathbf{z})}$ that was infected by some
$x_a \in\mathcal{C}_{f_\xi(\mathbf{z})}$. But then $w_{\mathbf
{z}}(x_a, x_b)>\xi$, and by {(ii)}, this is when they first
interact after $x_a$ is infected. Hence no infection could have been
passed. Thus no particle in $B$ is infected.

The following lemma generalises Lemma~\ref{ineractionapprox}

\begin{lemma}\label{ineractionapproxPsi}
Assume the conditions of Section~\ref{assumptionssection}.

Suppose $\mathbf{z}$ is good. Define a partition of the particles into
sets $A$ and $B$ where $A=V(\mathcal{C}_{f_\xi(\mathbf{z})})$. Let
$\mathcal{S}=E_{\mathbf{z}}(A) \cup E_{\mathbf{z}}(A:B)$. Then
\[
\Pr \bigl(E_{\Upsi'}(A)\cup E_{\Upsi'}(A:B)=\mathcal{S} \bigr) =
\bigl(1+o(1) \bigr)\Pr \bigl(E_{\Lambda}(A)\cup E_{\Lambda}(A:B)=
\mathcal{S} \bigr).
\]
\end{lemma}

\begin{pf}
The proof follows the same pattern as that of Lemma~\ref
{ineractionapprox}. Using Lemma~\ref{TimeDistPsi}, $E_{\Upsi'}(A)\cup
E_{\Upsi'}(A:B)=\mathcal{S}$ defines a sequence of interactions
$\sigma
(\mathbf{z})= ((e_{(1)},\tau_1), (e_{(2)}, \tau_2), \ldots,
(e_{ \bigl({|\mathcal{S}|\choose2} \bigr)}, \tau_{{|\mathcal
{S}|\choose2}}) )$ and active edges sets $A_0, A_1, \ldots,
A_{|\mathcal
{S}|}$. We get an equivalent to \eqref{faoifs},
%
\begin{equation}
\Pr \bigl(E_{\Upsi'}(A)\cup E_{\Upsi'}(A:B)=\mathcal{S} \bigr)=
\prod_{j=0}^{
{|\mathcal{S}|\choose2}-1}p\bigl(1-|A_j|p\bigr)^{\t_{j+1}-\t_j},
\label{l8kj5xdv}
\end{equation}
where
\[
p=\frac{1}{\theta_r n}{\biggl(1+O{\biggl(\frac{1}{n^{\Omega
(1)}}+\frac
{1}{(k \ln n)^{\Omega^+(1)}}
\biggr)}\biggr)}.
\]
We eventually get
\[
\Pr \bigl(E_{\Upsi'}(A)\cup E_{\Upsi'}(A:B)=\mathcal{S}
\bigr)=p^{{|\mathcal
{S}|\choose2}}\exp \biggl\{- \bigl(1+O \bigl(k^2p \bigr)
\bigr)p \sum_{e \in\mathcal
{S}}w_{\mathbf
{z}}(e) \biggr\}.
\]

Again, $q=\frac{\psi}{\theta_r n}$ where here $\psi=1$ because
$\Row=1$.
\begin{eqnarray*}
\Pr \bigl(E_{\Lambda}(A)\cup E_{\Lambda}(A:B)=\mathcal{S}
\bigr)&=&q^{|\mathcal
{S}|}(1-q)^{\sum_{e \in\mathcal{S}}w_{\mathbf{z}}(e)-|\mathcal
{S}|}
\\
&=&q^{|\mathcal{S}|}(1-q)^{-|\mathcal{S}|}\exp \biggl\{- \bigl(1+O(q) \bigr)q\sum
_{e
\in\mathcal{S}}w_{\mathbf{z}}(e) \biggr\}.
\end{eqnarray*}
Then taking the ratio and by the same reasoning as for Lemma~\ref
{ineractionapprox}, the result follows.
\end{pf}

An event $\mathcal{E}$ for $\mathcal{C}_{f_\xi(\Lambda)}$ (or
$\mathcal
{C}_{f_\xi(\Upsi)}$) is a subset of the possible realisations of
$\mathcal{C}_{f_\xi(\Lambda)}$ (resp., $\mathcal{C}_{f_\xi(\Upsi)}$).
When we speak of a realisation, we include edge weights as well as
graphical structure. For example, $S_g = \{\mathcal{C}_{f_\xi(\mathbf
{z})} \dvtx  \mathbf{z}\ \mathrm{good}\}$ is an event.

The following is a generalisation of Lemma~\ref{whpEqual} (proof in the
\hyperref[app]{Appendix}).

\begin{lemma}\label{probConv}
Assume the conditions in Section~\ref{assumptionssection}.

Suppose there is a sequence of events $(\mathcal{E}_n)_{n \geq1}$ and
a constant $p_c$ such that $\Pr(\mathcal{C}_{f_\xi(\Lambda)} \in
\mathcal{E}_n)\rightarrow p_c$ as $n \rightarrow\infty$. Then $\Pr
(\mathcal{C}_{f_\xi(\Upsi)} \in\mathcal{E}_n)\rightarrow p_c$ as $n
\rightarrow\infty$.
\end{lemma}

\begin{pf*}{Proof of Theorem~\ref{epidemic} when $\Row=1$}
Since $f_{\xi}(\Lambda)$ is an Erd\H{o}s--R\'enyi random graph
$\mathcal
{G}_{k,\hat{q}}$ where $\hat{q}=1-(1-\frac{1}{\theta_r n})^\xi$,
we can
apply standard results (see, e.g., \cite{Remco,Lukzak}) to
determine the size of $\mathcal{C}_{f_{\xi}(\Lambda)}$. We address
cases (i) and (iii) of Theorem~\ref{epidemic} first. (i) If $k\hat
{q}<1$, then there is a constant $\alpha$, such that $\Pr\brac
{|\mathcal
{C}_{f_\xi(\Lambda)}|/\ln k \leq\alpha} \rightarrow1$. Applying Lemma~\ref{probConv} for the events $\mathcal{E}_n=\{\mathcal{C}_{f_\xi
(\mathbf{z})} \dvtx |\mathcal{C}_{f_\xi(\mathbf{z})}|/\ln k \leq\alpha
\}$
gives the result. (iii) If $k\hat{q}>(1+\varepsilon)\ln k$ where
$\varepsilon$ is any positive constant then $\Pr(|\mathcal
{C}_{f_{\xi
}(\Lambda)}|=k) \rightarrow1$. Applying Lemma~\ref{probConv} for the
events $\mathcal{E}_n=\{\mathcal{C}_{f_\xi(\mathbf{z})} \dvtx |\mathcal
{C}_{f_\xi(\mathbf{z})}|=k\}$ gives the result.

Case (ii) requires slightly more consideration.
We have $\mathcal{G}_{k,\hat{q}}$ where $k\hat{q} \rightarrow c$ for
some constant $c>1$. Denote the largest component by $\mathcal{C}_1$
and let $\beta=\beta(c)$ denote the unique solution in $(0,1)$ of the
equation $\beta+e^{-\beta c}=1.$ Then every $\nu\in(\frac{1}{2},1)$
there exists a $\delta=\delta(\nu,c)$ such that
$
\Pr\brac{\llvert  \frac{|\mathcal{C}_1|}{k}-\beta\rrvert \leq
k^{-\nu
}}=1-O(k^\delta)
$; see, for example, \cite{Remco}. With high probability, all other
components have size $O(\ln k)$.
By the symmetry of $\Lambda$, $\Pr(\mathcal{C}_{f_\xi(\Lambda
)}=\mathcal
{C}_1)=\Pr(x_0 \in\mathcal{C}_1)=\E[|\mathcal{C}_1|]/k\rightarrow
\beta
$ as $k \rightarrow\infty$.
Hence,
\begin{eqnarray*}
&&\Pr{\biggl( \biggl\llvert \frac{|\mathcal{C}_{f_\xi(\Lambda)}|}{k}-\beta \biggr\rrvert \leq
k^{-\nu}\biggr)}
\\
&&\qquad = \Pr{\biggl( \biggl\llvert \frac{|\mathcal{C}_{1}|}{k}-\beta \biggr\rrvert \leq
k^{-\nu
} \wedge\mathcal{C}_{f_\xi(\Lambda)}=\mathcal{C}_1
\biggr)}\\
&&\qquad\quad{}+ \Pr{\biggl( \biggl\llvert \frac{|\mathcal{C}_{f_\xi(\Lambda)}|}{k}-\beta \biggr\rrvert \leq
k^{-\nu} \wedge\mathcal{C}_{f_\xi(\Lambda)}\neq\mathcal
{C}_1\biggr)}
\\
&&\qquad \rightarrow\beta.
\end{eqnarray*}
Applying Lemma~\ref{probConv} for the events $\mathcal{E}_n=\{
\mathcal
{C}_{f_\xi(\mathbf{z})} \dvtx |\mathcal{C}_{f_\xi(\mathbf{z})}/k-\beta
|\leq
k^{-\nu}\}$, we get
$\Pr\brac{\llvert  \frac{|\mathcal{C}_{f_\xi(\Upsi)}|}{k}-\beta
\rrvert \leq k^{-\nu}} \rightarrow\beta$, completing the proof of this case.
\end{pf*}


\section{\texorpdfstring{Extending to the general case: SI(R) with $\rho\leq
1$}{Extending to the general case: SI(R) with $rho<=1$}}\label
{TheGeneralCase}

This section deals with general $\xi$ and $\Row\leq1$.

\subsection{\texorpdfstring{A heuristic treatment of a two-particle system.}
{A heuristic treatment of a two-particle system}}\label
{twoparticles}
We introduce the approach we shall use by giving an informal treatment
for a system with only two particles. In subsequent sections, we will
give a more rigorous analysis for $k$ particles, where, as per the
assumptions in Section~\ref{assumptionssection}, $k \leq n^\epsilon
$, for $\epsilon
$ a sufficiently small constant.

Let $x$ and $y$ be the two particles, with $x$ being the initial
infective and $y$ being susceptible. We allow $\xi<\infty$ and/or
$\Row
< 1$. The former conditions means that $y$ may never get infected, the
latter condition means that it may take more than one meeting between
$x$ and $y$ before an interaction takes place. Note that if $x$ and $y$
were at the same vertex at time $t$, and happen to move to the same
neighbouring vertex in the next step, then this counts as another
meeting, with another coin flip to determine if an interaction takes place.

Now, suppose $x$ and $y$ have just stepped to the same vertex $v$. With
probability $\Row$ there will be an interaction. After this, they will
move again, either to the same neighbour of $v$ with probability $1/r$
or to different neighbours with probability $(r-1)/r$. Let
\begin{eqnarray*}
\phi&=& \Pr(\mbox{No $xy$ interaction before they move apart})
\\
&=& \sum_{i\geq1} (1-\Row )^i \biggl(
\frac{1}{r} \biggr)^{i-1} \biggl(1-\frac{1}{r} \biggr) =
\frac{(1-\Row)(r-1)}{r-1+\Row}.
\end{eqnarray*}

Recall from Section~\ref{TypicalGraphsSection} that a vertex $v$ is
\emph{treelike} if there
is no cycle in the subgraph $G[v,L_1]$ induced by the set of
vertices within (graph) distance $L_1=\rdown{\epsilon _1\log_rn}$ of $v$,
where $\epsilon _1>0$ is a sufficiently small constant. The following
lemma is
from \cite{CFR-Mult}:

\begin{lemma}[(\cite{CFR-Mult})]\label{2-meet}
Let $G$ be a typical $r$-regular graph, and let $v$ be a vertex of~$G$, treelike to depth $L_1=\rdown{\epsilon _1\log_r n}$. Suppose
that at
time zero, two independent random walks $(\cW_1,\cW_2)$ start from
$v$. Let $(a,b)$ denote the position of the particles at any step.
Let $S=\set{(u,u)\dvtx u \in V}$. Let $f$ be the probability of a first
return to $S$ within $T=O(k \ln n)$ steps given that the walks leave $v$
by different edges at time zero. Then
\[
f = \frac{1}{(r-1)^2}+O \bigl(n^{-\OL} \bigr).
\]
\end{lemma}

Using this lemma, let
%
\begin{eqnarray}\label{phideriv}
\phi_T &=& \Pr(\mbox{No $xy$ interaction before being apart more
than $T$ time steps})\hspace*{-27pt}
\nonumber
\\[-4pt]
\\[-12pt]
\nonumber
&=& \sum_{i \geq1}\phi^if^{i-1}
(1-f ) = \frac{\phi
(1-f )}{1-\phi f}.\hspace*{-27pt}
\end{eqnarray}

If two sequences $A_n, B_n$ are such that $A_n/B_n \rightarrow1$ as $n
\rightarrow\infty$, we write $A_n \sim B_n$.
The rest of this section will be implicitly justified in the detailed
treatment in Section~\ref{SIrholeq1case}.

Recall $\psi$ was defined in \eqref{defpsi}, and observe that $\Row
\leq
\psi\leq1$ with $\psi=1$ if and only if $\Row=1$.
Now, assuming $x$ and $y$ start at the same vertex, then as will be
seen from Lemma~\ref{2-interact},
%
\begin{eqnarray}\label{psideriv}
&&\Pr(\mbox{$xy$ interaction occurs within $T$ time steps})
\nonumber
\\
&&\qquad\sim\Pr(\mbox{$xy$ interaction occurs before $x$ and $y$ have been}
\nonumber\\[-8pt]\\[-8pt]\nonumber
&&\hspace*{160pt}\mbox{apart more than $T$ steps})
\\
&&\qquad =1-\phi_T = 1-\frac{\phi (1-f )}{1-\phi f} \sim\frac
{\Row(r-1)}{r-2+\Row} = \psi.\nonumber
\end{eqnarray}

\begin{definition}[({$\ell$-distinct})]\label{TdistinctDef}
A sequence $(t_1, t_2, \ldots)$ is \emph{$\ell$-distinct} if $t_1
\geq
\ell$ and $t_{i+1}-t_i \geq\ell$.
\end{definition}

Clearly, there can be at most $t/\ell$ $\ell$-distinct meetings in
$[0,t]$, and assuming $i \leq t/\ell$,
%
\begin{equation}
\Pr \bigl(\mbox{there are $i$ $\ell$-distinct meetings in $[0,t]$} \bigr) \sim
\pmatrix{t
\cr
i}p^{i}(1-p)^{t-i},
\end{equation}
where $p=(1+o(1))\frac{1}{\theta_r n}$ is from \eqref{longformp}.
Hence, it is seen to be approximately distributed as $\operatorname{Binom}(t,\frac{1}{\theta_r n})$. The probability that there are no
interactions in any of the $i$ intervals $[t_j, t_j+T]$ where $t_j$ is
the time of the $j$th $\ell$-distinct meeting is $(1-\psi)^i$. Thus
\begin{eqnarray*}
&& \Pr \bigl(\mbox{there are no interactions in the period $[0,t]$} \bigr)
\\
&&\qquad \sim \sum_{i =0}^{t/\ell}\pmatrix{t
\cr
i}
\bigl(p(1-\psi ) \bigr)^{i}(1-p)^{t-i}
\\
&&\qquad \sim (1-\psi p)^{t}.
\end{eqnarray*}
Hence
%
\begin{equation}
\Pr(\mbox{$y$ gets infected within time $\xi$}) \sim 1-(1-\psi
p)^{\xi}.\label{janf}
\end{equation}
When $\Row=\psi=1$, \eqref{janf} looks similar to the bracketed terms
in \eqref{np}. This is, of course, not a coincidence since the
bracketed term in \eqref{np} is essentially the probability that an
infection is passed between a pair of particles if one of them had been
infected. Therefore, $\Phi$ is effectively the expected number of other
particles that are infected by a particular particle.


\subsection{\texorpdfstring{Allowing $\rho<1$.}{Allowing $rho<1$}}\label{SIrholeq1case}
In this section, we formalise some of the ideas of Section~\ref{twoparticles}, extended to $k$ particles.

Let us first redefine $\ell$ to be $\ell= 2(T+T^3)+T$ where $T$ is a
maximal mixing time.

Consider a period $[0,t]$. Let $\ull\t=(t_1,\ldots,t_i, t_{i+1})$ be a
sequence such that (i) $\sum_{s=1}^{i+1}t_s=t$, (ii) $t_s \geq\ell$ for
each $s \in\{1,2,\ldots,i\}$, (iii) $t_i+t_{i+1} \geq\ell+T$ and (iv)
$0 \leq t_{i+1} \leq T$.

Let $t^*_0=0$ and $t^*_s=t^*_{s-1}+t_s$ for $s=1,2,\ldots,i$. We use
$\ull\t$ to represent having $i$ $\ell$-distinct meetings in the period
$[0,t]$ with a first interaction at step $t$ (recall definition of
$T$-distinct given in Definition~\ref{TdistinctDef}). Specifically, let
$A$ be a set of active particle pairs and let $(x,y) \in A$. Denote by
$\mathcal{C}_{(x,y)}(\ull\t)$ the following event: (1) at each time
$t^*_s$ a single particle pair $(a,b)_s \in A$ meets and does not
interact in the period $[t^*_s, t^*_s+T]$. (2) At time $t^*_i$, $(x,y)$
meet and interact at some point in the period $[t^*_i, t]$. (3) No
particle pair $(a,b) \in A$ meets in the periods $[t^*_s+\ell,
t^*_{s+1}-1]$, $s =0, 1, \ldots, i$.

Note, in (2), we do not specify which particle pair $(a,b)$ meet, only
that a single pair meet. Thus, there are $|A|^{i-1}$ possible sequences
of particles that satisfy this condition.

Let ${\ull\t}_i(t)$ be the set of all sequences $\ull\t=(t_1,\ldots,t_i,
t_{i+1})$ satisfying the above conditions, and let ${\ull\t}(t)=\bigcup_{i}{\ull\t}_i(t)$.

Define the event $\mathcal{C}_{(x,y)}(t)=\bigcup_{{\ull\t} \in{\ull\t
}(t)}\mathcal{C}_{(x,y)}(\ull\t)$. Thus $\mathcal{C}_{(x,y)}(t)$
represents, subject to conditions (i)--(iv), having some number of
noninteractive meetings of active particle pairs, and the first
interaction taking place between $x$ and $y$ at time $t$.

The following lemma is core to this section. The proof is rather long,
so we delay it till later.

\begin{lemma}\label{mainPart3Lemma}
Let $G$ be typical, and let $k \leq n^\epsilon $ for sufficiently
small $\epsilon $.

Suppose $\sqrt{n}\leq t \leq k^2 n \ln n$. For a set of active particle
pairs $A$ and $(x,y) \in A$,
\[
\Pr \bigl(\mathcal{C}_{(x,y)}(t) \bigr)=\psi p(1-|A|\psi
p)^t,
\]
where
%
\begin{equation}
p=\frac{1}{\theta_r n}{\biggl(1+O{\biggl(\frac{1}{n^{\Omega
(1)}}+\frac
{1}{(k \ln n)^{\Omega^+(1)}}
\biggr)}\biggr)}
\end{equation}
and
\[
\psi=\frac{\Row(r-1)}{r-2+\Row}.
\]
\end{lemma}

In the proofs of Lemmas \ref{ineractionapprox} and~\ref
{ineractionapproxPsi}, we dealt with epochs $\t_j$ and sets $A_j$ of
active edges. We assumed $\mathbf{z}$ was good and calculated the
probability of a sequence $\sigma(\mathbf{z})$. This sequence specified
that no particle pair in $A_j$ would interact in the period $[\tau_j,
\tau_{j+1}-1]$ then a particular pair $(x,y)$ interact at time $\tau
_{j+1}$. Because $\Row=1$, interactions coincided with meetings, and so
the sequence $\sigma(\mathbf{z})$ specified meetings and nonmeetings.
When $\Row<1$, there may be some number of noninteractive meetings of
active particle pairs in the period $[\tau_j, \tau_{j+1}-1]$ before the
$xy$ interaction finally takes place at time $\t_{j+1}$, and we will
need to take these into account when calculating the probability of
$\sigma(\mathbf{z})$.

Recall that the raison d'\^etre of \emph{good} $\mathbf{z}$, rather
than any arbitrary $\mathbf{z}$, was that having meetings take place
$\ell$ steps apart (i.e., being $\ell$-distinct) allowed us to apply
Lemma~\ref{ProbBxyLemma}. Although this lemma could tell us what
happens in the period $[\t_j+\ell, \t_{j+1}]$, we could not account for
what happens in the first $\ell$ steps after the epoch~$\t_j$. We
therefore modified the process in Section~\ref{UpsilonPrimeSection} so
that $\Row(t)=0$ in these periods we could not account for, thereby
guaranteeing that no interaction occurred and therefore, the
probability calculated by Lemma~\ref{ProbBxyLemma} was faithful to the
modified process. Hence we calculated probabilities for $\Upsi'$ rather
than $\Upsi$, and we then related the two by proving they are the same
w.h.p. Since Lemma~\ref{ProbBxyLemma} is the main tool in the proof of
Lemma~\ref{mainPart3Lemma}, we must do something similar.

In this section, we need to have a slightly different version of the
interaction graph $\Upsi'$, which we shall denote by $\Upsi^*$. To cope
with sequences $\ull\t$ satisfying conditions (i)--(iv), we modify the
process as follows; cf. Section~\ref{UpsilonPrimeSection}:

\begin{longlist}[(ii)]
\item[(i)] Set $\Row(t)=0$ for $t \in[1, \ell]$.
\item[(ii)] If at time $\t$ there was a meeting between one
active pair $(x,y)\in A$, do the following: (a) For $(x,y)$, set $\Row
(t)=0$ for $t \in[\t+T, \t+\ell]$. (b) For every other $(a,b) \in A$,
set $\Row(t)=0$ for $t \in[\t+1, \t+\ell]$.
\end{longlist}
Note, there remains the possibility of more than one active particle
pair meeting at a particular time step. With high probability, this
will not happen, as stated in Lemma~\ref{GoodLemExt}. However, for
completeness, we will include the following component to the above
modification, which represents a ``failure'':
\begin{longlist}[(iii)]
\item[(iii)] If at time $\t$ there were meetings of more than
one active pair in $A$, then set $\Upsi^*=(-1,-1,\ldots,-1)$ and
terminate the construction.
\end{longlist}
Note, we did not need {(iii)} in the $\Row=1$ analysis because
we calculated meetings based on $\mathbf{z}$ being good, which required
that no two pairs of active particles simultaneously meet, in
accordance with part {(a)} of the definition of good (Definition~\ref{goodDef}). In this section, we consider meetings where no
interaction takes places, and modification {(iii)} serves as a
technical convenience.

The implication of these modifications is that in the application of
Lem\-ma~\ref{ProbBxyLemma} in the proof of Lem\-ma~\ref{mainPart3Lemma},
there were no interactions in the periods that we could not account for.

$\Upsi^*$ is essentially a generalisation of $\Upsi'$. For $\Upsi'$, we
did not need to stipulate {(ii)(a)} above, since if $(x,y)$ met,
they interacted, and further interaction between them thereafter had no
bearing on edge weight $w_{\Upsi'}(x,y)$. Hence {(ii)(a)} would
have been redundant.

The following lemma allows us to say that under these new
modifications, w.h.p., $\Upsi=\Upsi^*$. The proof is in the \hyperref[app]{Appendix}.

\begin{lemma}\label{GoodLemExt}
Let $G$ be typical, and let $k \leq n^\epsilon$ for sufficiently small
$\epsilon$.

With high probability:
\begin{longlist}[(a)]
\item[(a)] only one pair meet at a time, and no other pair meet within
$\ell= 3T +2T^3$ steps;
\item[(b)] any active pair that meets at some step $\t$ and does not
interact in the period $[\t,\t+T-1]$ does not meet in the period $[\t
+T, \t+\ell]$.
\end{longlist}
\end{lemma}

\begin{corollary}\label{UpsiesEqualMostGeneral}
With high probability, $\Upsi=\Upsi^*$.
\end{corollary}

Consequently, if for a sequence of events $\mathcal{E}_k$ we determine
that in $\Upsi^*$, $\Pr(\mathcal{E}_k)\rightarrow p_c$ where $p_c$
is a
constant, then it will also be the case in $\Upsi$ that $\Pr(\mathcal
{E}_k)\rightarrow p_c$.

Observe an assumption in Lemma~\ref{mainPart3Lemma} is $\sqrt{n} \leq t
\leq k^2 n \ln n$. The RHS inequality follows from the definition of
good (Section~\ref{interapprox}). The LHS inequality is an extra
condition we must impose; thus we redefine Definition~\ref{goodDef}
part {(a)} as follows: \emph{If $\Upsi=\mathbf{z}$, then
none of
the ${k\choose2}$ interactions that form the edges of $\Upsi$ took
place within $\sqrt{n}$ steps of each other}.
Showing that Lemma~\ref{UpsiesEqual} still holds is straightforward: by
Lemma~\ref{L3} the probability of any of the at most ${k\choose2}$
particles meeting within $\sqrt{n}$ steps is $O(k^2\sqrt {n}/n)=O(k^2/\sqrt{n})$. Taken over all [at most ${k\choose2}$]
interactions, this is $O(k^4/\sqrt{n})=o(1)$ when $k \leq n^{\epsilon}$
for small enough $\epsilon$.

We remind the definition of $\Lambda$ given in Section~\ref{interapprox} stipulates i.i.d. edge weights distributed as $\Ge(q)$
where $q=\psi/(\theta_r n)$.

\begin{lemma}\label{interactionPart3}
Assume the conditions of Section~\ref{assumptionssection}.

Suppose $\mathbf{z}$ is good. Define a partition of the particles into
sets $A$ and $B$ where $A=V(\mathcal{C}_{f_\xi(\mathbf{z})})$. Let
$\mathcal{S}=E_{\mathbf{z}}(A) \cup E_{\mathbf{z}}(A:B)$. Then
\[
\Pr \bigl(E_{\Upsi^*}(A)\cup E_{\Upsi^*}(A:B)=\mathcal{S} \bigr) =
\bigl(1+o(1) \bigr)\Pr \bigl(E_{\Lambda}(A)\cup E_{\Lambda}(A:B)=
\mathcal{S} \bigr).
\]
\end{lemma}

\begin{pf}
Observe that Lemma~\ref{TimeDistPsi} holds for $\Upsi^*$ in the same
way as it does for $\Upsi'$. It is stated and proved in terms of
interactions; noninteractive meetings are irrelevant. Hence, we can
now follow the pattern of the proof of Lemma~\ref{ineractionapproxPsi}:
$E_{\Upsi^*}(A)\cup E_{\Upsi^*}(A:B)=\mathcal{S}$ defines a sequence of
interactions $\sigma(\mathbf{z})= ((e_{(1)},\tau_1), (e_{(2)},
\tau
_2), \ldots, (e_{ \bigl({|\mathcal{S}|\choose2} \bigr)}, \tau
_{{|\mathcal{S}|\choose2}}) )$ and active edges sets $A_0, A_1,
\ldots,
A_{|\mathcal{S}|}$.

Consider the period between $\t_j$ and $\t_{j+1}$ for $1 \leq j \leq
{|\mathcal{S}|\choose2}-1$. By construction of $\Upsi^*$,
specifically, by the process modifications {(i)} and
{(ii)} above, meetings in the period $[\t_j, \t_{j+1}]$ must be of the
form $\ull\t$ satisfying conditions (i)--(iv) above, except those for
which $t_i+t_{i+1} < \ell+T$, that is, those for which $t^*_{i-1} \in
[t-(\ell+T), t-\ell]$. This window of size $T$ can be ignored since,
from the stationary distribution, the processes has probability
$O(k^2T/n)$ of falling into it. Taken over all [at most ${k\choose2}$]
periods this is still $o(1)$ for $k \leq n^\epsilon $ and $\epsilon $
small enough.

We apply Lemma~\ref{mainPart3Lemma} to each period $[\t_j, \t_{j+1}]$
get to get an equivalent to \eqref{l8kj5xdv}:
%
\begin{equation}
\Pr \bigl(E_{\Upsi^*}(A)\cup E_{\Upsi^*}(A:B)=\mathcal{S} \bigr)=
\prod_{j=0}^{
{|\mathcal{S}|\choose2}-1}\psi p\bigl(1-|A_j|
\psi p\bigr)^{\t_{j+1}-\t_j},
\end{equation}
where
\[
p=\frac{1}{\theta_r n}{\biggl(1+O{\biggl(\frac{1}{n^{\Omega
(1)}}+\frac
{1}{(k \ln n)^{\Omega^+(1)}}
\biggr)}\biggr)}.
\]

Recalling $q=\frac{\psi}{\theta_r n}$, the rest of the proof follows in
the same way as in Lemma~\ref{ineractionapproxPsi} (and in turn, Lemma~\ref{ineractionapprox}).
\end{pf}

Note, without modification {(iii)}, which sets $\Upsi^*$ to a
vector of $-1$'s we would have to consider sequences which give rise to
$\Upsi^*=\mathbf{z}$ and in which more than one active particle meet.
Thus modification {(iii)} is a technical convenience.

With the proof of Lemma~\ref{interactionPart3}, the proofs of Theorems
\ref{epidemic} and \ref{completiontimetheorem} follow for the general
case $\Row\leq1$ as they did for the special case $\Row=1$. We address
Theorem~\ref{completiontimetheorem} first.

\begin{pf*}{Proof of Theorem~\ref{completiontimetheorem}}
Recall that Lemma~\ref{UpsiesEqual} is general, holding for $\xi\leq
\infty$ and $\Row\leq1$. Thus, if $\xi=\infty$, $\Upsi$ remains good,
w.h.p., when $\Row<1$. By Corollary~\ref{UpsiesEqualMostGeneral}, the
same is true for $\Upsi^*$. Consequently, by Lemma~\ref
{interactionPart3}, $\Lambda$ is good w.h.p. Thus, we have all three
components of Corollary~\ref{UpsiesAREequal} holding for $\Row\leq1$.
As such, Lemma~\ref{whpEqual} holds for $\xi=\infty$ and $\Row\leq1$.
Therefore, the proof of Theorem~\ref{completiontimetheorem} given in
Section~\ref{interapprox} holds with $\lambda=\psi/(\theta_r n)$.
\end{pf*}

\begin{pf*}{Proof of Theorem~\ref{epidemic}}
Observe that the proof of Lemma~\ref{probConv} used Corollary~\ref
{UpsiesAREequal} and Lemma~\ref{ineractionapproxPsi}. We have
generalised both of these in this section, so Lemma~\ref{probConv}
holds for the general case $\xi\leq\infty$, $\Row\leq1$ in the same
way, with the edges of $\Lambda$ being i.i.d. as $\Ge(\psi/(\theta_r
n))$. Consequently, the general case of Theorem~\ref{epidemic} is
justified in the same was as the special case was in Section~\ref{SIRRho1} [observe that now $f_{\xi}(\Lambda)$ is an Erd\H{o}s--R\'enyi
random graph $\mathcal{G}_{k,\hat{q}}$ where $\hat{q}=1-(1-\frac
{\psi
}{\theta_r n})^\xi$].
\end{pf*}

Before we proceed to prove Lemma~\ref{mainPart3Lemma}, we require the
following:

\begin{lemma}\label{2-interact}
Let $G$ be typical and let $k \leq n^\epsilon$ for sufficiently small
$\epsilon$.

Suppose that at time zero, two particles $x,y$ positioned on a treelike
vertex $v$ interact with probability $\Row$. Let $\psi'$ be the
probability of an $xy$ interaction in $[0,T-1]$ where $T$ is a mixing
time. Then
%
\begin{equation}
\psi'=\psi \bigl(1- O \bigl(n^{-\Omega(1)} \bigr) \bigr),
\label{psi'def}
\end{equation}
where $\psi=\frac{\Row(r-1)}{r-2+\Row}$.
\end{lemma}

\begin{pf}
$x$ and $y$ start at the same vertex at time $t=0$, and interact with
probability $\Row$ at time $0$. Suppose the first $xy$ interaction
occurs at time $\t\geq0$. Now suppose $x$ and $y$ are incident (i.e.,
at the same vertex) at times $0=t_0, t_1, \ldots.$ Let $t_{M}$ be the
smallest $t_r$ in this sequence such that $t_{M+1}-t_{M}\geq T$. Hence,
$t_{M}+T-1$ is the first time that $x$ and $y$ have been apart $T-1$ steps.

We first demonstrate that $\Pr(\t<T) = \Pr(\t< t_{M}+T) +o(1)$. Observe
\[
\Pr(\t\geq t_{M}+T) \leq\Pr(\t\geq T) \leq\Pr(\t\geq
t_{M}+T)+ \Pr \bigl( \bigl\{x, y\mbox{ meet in }[T,2T] \bigr\}
\bigr).
\]
Using \eqref{atv},
\begin{eqnarray*}
\Pr \bigl(x, y\mbox{ meet in }[T,2T] \bigr)&=&1-\frac{(1+O(T\p_\gamma
))}{
{(1+(1+O(T\p_\gamma ))\p_\gamma /R_\gamma) }^{2T}}-O
\bigl(T^2\pi {_\gamma }e^{-\lambda T} \bigr)
\\
&=& O \bigl(T^2\pi_\gamma \bigr)
\\
&=& O \bigl((k\ln n)^2/n \bigr).
\end{eqnarray*}

Thus
\[
\Pr(\t\geq t_{M}+T) \leq\Pr(\t\geq T) \leq\Pr(\t\geq
t_{M}+T)+O \bigl((k\ln n)^2/n \bigr),
\]
that is, $\Pr(\t< t_{M}+T)=\Pr(\t<T)+O((k\ln n)^2/n)$.

Recall Section~\ref{twoparticles}. To get the correcting factor in
\eqref{phideriv}, observe that the sum $\sum_{i \geq1}\phi
^if^{i-1} (1-f )$ assumes every vertex at which the particles
part is tree-like. Let $\cW(t)$ be the walks on $G$, and let $\cX(t)$
be a walk on an infinite $r$-regular tree $\cT$ rooted at the start
vertex $v$, which is assumed to be tree-like. We couple $\cW$ and $\cX$
until time $L_1$. Since $G$ and $\cT$ have the same structure out to
$L_1$, the two processes are identical until $t=L_1$. Let
$Y_t=\operatorname
{dist}(x,y)$ in $G$. It is shown in \cite{CFR-Mult}, proof of Lemma~17, that $\Pr(Y_{L_1}\leq L_1/2) = O(n^{-\Omega^+(1)})$ where the
$\Omega^+(1)$ is an arbitrarily large constant. It is also shown that,
subject to $k\leq n^\epsilon $ for a sufficiently small $\epsilon $,
$\Pr(\mbox{the
walks meet in $[L_1, T]$ and $Y_{L_1}>L_1/2$})=O(n^{-\Omega(1)})$. Thus
\[
\phi_T=\frac{\phi (1-f )}{1-\phi f}+O \bigl(n^{-\Omega(1)} \bigr).
\]

Now we include the error term for the asymptotic equality in line
\eqref{psideriv}.
%
\begin{eqnarray}\label{f8sdg8s46d8gs}
1-\phi_T&=&\frac{\Row(r-1)}{r-2+\Row} \bigl(1+O \bigl(n^{-\Omega
(1)}
\bigr) \bigr)
\nonumber
\\[-8pt]
\\[-8pt]
\nonumber
&=& \psi \bigl(1+O \bigl(n^{-\Omega(1)} \bigr) \bigr)=\Pr(\t<
t_{M}+T).
\end{eqnarray}
Note in \eqref{f8sdg8s46d8gs} we have used the assumption that $\Row$
is a constant to absorb (functions of) it into the $O$ term. Thus
\begin{eqnarray*}
\Pr(\t< T)&=&\Pr(\t< t_{M}+T)-O \biggl(\frac{(k\ln n)^2}{n} \biggr)\\
&=&\psi
\bigl(1+O \bigl(n^{-\Omega(1)} \bigr) \bigr)-O \biggl(\frac{(k\ln
n)^2}{n}
\biggr).
\end{eqnarray*}
$\Row$ is a constant, which means $\psi$ is a constant, and so defining
$\psi'=\Pr(\t< T)$,
\[
\psi'=\psi \biggl(1- O \biggl(\frac{1}{n^{\Omega(1)}} \biggr) \biggr).
\]
\upqed\end{pf}

\begin{pf*}{Proof of Lemma~\ref{mainPart3Lemma}}
Writing
\[
\Pr \bigl(\mathcal{C}_{(x,y)} \bigl(\ull\t_i(t) \bigr)
\bigr)= \sum_{\ull\t\in\ull\t
_i(t)}\Pr \bigl(\mathcal{C}_{(x,y)}(
\ull \t) \bigr),
\]
we have
%
\begin{eqnarray}\label{g46dzxg}
\Pr \bigl(\mathcal{C}_{(x,y)}(t) \bigr)&=&\sum
_{i \geq1}\Pr \bigl(\mathcal {C}_{(x,y)} \bigl(\ull
\t_i(t) \bigr) \bigr)
\nonumber
\\[-8pt]
\\[-8pt]
\nonumber
&=&\sum_{i=1}^{(k \ln n)^5}\Pr \bigl(
\mathcal{C}_{(x,y)} \bigl(\ull\t_i(t) \bigr) \bigr) + \sum
_{i \geq(k \ln n)^5+1}\Pr \bigl(\mathcal{C}_{(x,y)} \bigl(
\ull \t_i(t) \bigr) \bigr).
\end{eqnarray}

We shall focus on the first sum, returning to the second later.

We can write $\psi'=\sum_{i=0}^{T-1}\Row_i$ where $\Row_i$ is the
probability that the first interaction happens at step $i$. So, for
example, $\Row_0=\Row$. For a given $\ull\t$, let $\Row_{\ull\t
}=\Row
_{t_{i+1}}$.

We shall calculate $\Pr(\mathcal{C}_{(x,y)}(\ull\t))$. By Lemma~\ref
{ProbBxyLemma}, the probability that no pair in $A$ meet in the period
$[t^*_{s-1}+\ell,t^*_s-1]$, then some \emph{particular} pair $(a,b)
\in
A$ meet (and no others do) at time $t^*_s$ is given by
\[
\Pr \bigl(\cB_{(a,b)}(t_s) \bigr)={\biggl(1+O{\biggl(
\frac{1}{n^{\Omega
(1)}}+\frac
{1}{(k \ln n)^{\Omega^+(1)}}+\frac{k^4t_s}{n^2}\biggr)}
\biggr)}p\bigl(1-|A|p\bigr)^{t_s},
\]
where
\[
p=\frac{1}{\theta_r n}{\biggl(1+O{\biggl(\frac{1}{n^{\Omega
(1)}}+\frac
{1}{(k \ln n)^{\Omega^+(1)}}
\biggr)}\biggr)}.
\]
Since $t_s \leq t \leq k^2 n \ln n$, we can write $\Pr(\cB
_{(a,b)}(t_s))=p(1-|A|p)^{t_s}$.
Taking the product over all $i$ $\ell$-distinct meetings, we get
%
\begin{equation}
\prod_{s=1}^{i}p\bigl(1-|A|p\bigr)^{t_s}.\label{fodjosi454}
\end{equation}
In the first $i-1$ of these the pair do not interact in the following
$T$ steps, so we multiply \eqref{fodjosi454} by $(1-\psi')^{i-1}$, and
the interaction happens at time $t$, which is $t_{i+1} \leq T$ after
$t^*_1$, so we multiply by $\Row_{t_{i+1}}$. Thus we get
%
\begin{equation}
\Row_{\ull\t} \bigl(1-\psi' \bigr)^{i-1}\prod
_{s=1}^{i}p\bigl(1-|A|p\bigr)^{t_s}.\label{poskd}
\end{equation}

In the above expression, a particular pair is specified at each $\ell
$-distinct meeting. We wish to fix only the final pair, giving a total
of $|A|^{i-1}$ possible cases. Each case has the above probability, so
in total, we get
\begin{eqnarray*}
\Pr \bigl(\mathcal{C}_{(x,y)}(\ull\t) \bigr)&=&|A|^{i-1}
\Row_{\ull\t} \bigl(1-\psi ' \bigr)^{i-1}\prod
_{s=1}^{i}p\bigl(1-|A|p\bigr)^{t_s}
\\
&=&p \bigl( \bigl(1-\psi' \bigr)|A|p \bigr)^{i-1}\bigl(1-|A|p\bigr)^{t}
\Row_{\ull\t}\bigl(1-|A|p\bigr)^{-t_{i+1}}.
\end{eqnarray*}

Since $(1-|A|p)^{-t_{i+1}}=1+O(Tk^2/n)$, it can be absorbed into the
correcting factor in $p$, so we have
\[
\Pr \bigl(\mathcal{C}_{(x,y)}(\ull\t) \bigr)=\Row_{\ull\t}p
\bigl( \bigl(1-\psi ' \bigr)|A|p \bigr)^{i-1}\bigl(1-|A|p\bigr)^{t}.
\]

Observe that for a given $t$, fixing $t_1, \ldots, t_{i-1}$ and allowing
$t_{i+1}$ to vary from $0$ to $T$, determines $t_i$. Letting
$t^{(r)}_i=t-(t_1+\cdots+t_{i-1}+r)$, we have
\begin{eqnarray*}
\sum_{r=0}^{T}\Pr \bigl({\ull\t}=
\bigl(t_1, \ldots, t_{i-1}, t^{(r)}_i,
r \bigr) \bigr)&=&p \bigl( \bigl(1-\psi' \bigr)|A|p
\bigr)^{i-1}\bigl(1-|A|p\bigr)^{t} \sum_{r=0}^T
\Row_r
\\
&=&\psi'p \bigl( \bigl(1-\psi' \bigr)|A|p
\bigr)^{i-1}\bigl(1-|A|p\bigr)^{t}.
\end{eqnarray*}

Recall ${\ull\t}_i(t)$ is the set of all sequences $\ull\t
=(t_1,\ldots,t_i, t_{i+1})$ satisfying conditions (i)--(iv). Let ${\ull\t
}^{\star}_i(t) \subset{\ull\t}_i(t)$ be those for which $t_{i+1}=0$.
%
\begin{eqnarray}\label{g5sdg46s5}
\Pr \bigl(\mathcal{C}_{(x,y)} \bigl(\ull\t_i(t) \bigr)
\bigr) &=&\sum_{\ull\t\in{\ull\t}_i(t)} \Pr \bigl(\mathcal{C}_{(x,y)}(
\ull\t) \bigr)
\nonumber
\\
&=&\sum_{\ull\t\in{\ull\t}^{\star}_i(t)}\psi'p \bigl( \bigl(1-
\psi ' \bigr)|A|p \bigr)^{i-1}\bigl(1-|A|p\bigr)^{t}
\\
&=& \bigl|{\ull\t}^{\star}_i(t)\bigr|\psi'p \bigl(
\bigl(1-\psi ' \bigr)|A|p \bigr)^{i-1}\bigl(1-|A|p\bigr)^{t}.\nonumber
\end{eqnarray}

Let $\mathcal{S}_{i}$ be the set of sequences for which $t_i<\ell+T$,
or in which there is some $t_s<\ell$ for $1\leq s \leq i-1$.
\begin{eqnarray*}
\bigl|{\ull\t}^{\star}_i(t)\bigr|&=&\pmatrix{t
\cr
i-1}-|
\mathcal{S}_{i}|
\\
&=& \pmatrix{t
\cr
i-1}-O{\biggl(\pmatrix{t
\cr
i-2}i\ell\biggr)}
\\
&=&\pmatrix{t
\cr
i-1} {\biggl(1-O\biggl(\frac{i^2\ell} {t-i}\biggr)\biggr)}
\\
&=&\pmatrix{t
\cr
i-1} {\biggl(1-O\biggl(\frac{(k \ln n)^{13}} {t}\biggr)\biggr)}.
\end{eqnarray*}

The above holds because we are only considering $i \leq(k \ln n)^5$,
and by assumption, $t \geq\sqrt{n}$, so when $k \leq n^\epsilon$ for
small enough $\epsilon$, we have $i =o(t)$. On the same basis, we can
further reduce the fraction to $1/n^{\Omega(1)}$, which can be absorbed
into the correcting factor of $p$.

Thus, continuing from \eqref{g5sdg46s5},
\begin{eqnarray*}
&&\sum_{i=1}^{(k \ln n)^5}\Pr \bigl(
\mathcal{C}_{(x,y)} \bigl(\ull\t_i(t) \bigr) \bigr)
\\
&&\qquad = \psi' p\sum_{i=1}^{(k \ln n)^5}\bigl(1-|A|p\bigr)^{i-1}
\pmatrix {t
\cr
i-1} \bigl( \bigl(1-\psi' \bigr)|A|p
\bigr)^{i-1}\bigl(1-|A|p\bigr)^{t-i+1}
\\
&&\qquad =\psi' p \sum_{i=1}^{(k \ln n)^5}
\pmatrix{t
\cr
i-1} \bigl( \bigl(1-\psi ' \bigr)|A|p
\bigr)^{i-1}\bigl(1-|A|p\bigr)^{t-i+1},
\end{eqnarray*}
since $(1-|A|p)^{i-1}=1+O(1/n^{\Omega(1)})$ when $i \leq(k \ln n)^5$.

Now
\begin{eqnarray*}
&&\sum_{i=1}^{t}\pmatrix{t
\cr
i-1} \bigl(
\bigl(1-\psi' \bigr)|A|p \bigr)^{i-1}\bigl(1-|A|p\bigr)^{t-i+1}\\
&&\qquad=
\bigl(1-\psi'|A|p \bigr)^t - \bigl( \bigl(1-
\psi' \bigr)|A|p \bigr)^t.
\end{eqnarray*}

Since $\sqrt{n} \leq t$, $((1-\psi')|A|p)^t= O((k^2/n)^{\sqrt{n}})
=O(1/n^{\Omega(\sqrt{n})})$.

Hence, putting the above into \eqref{g46dzxg}, we have
\begin{eqnarray*}
&&\Pr \bigl(\mathcal{C}_{(x,y)}(t) \bigr)\\
&&\qquad=\psi'p \bigl(1-
\psi'|A|p \bigr)^t - O\biggl(\frac {1} {n^{\Omega(\sqrt{n})}}\biggr)
\\
& &\qquad\quad{}+ \sum_{i \geq(k \ln n)^5+1}\Pr \bigl(\mathcal{C}_{(x,y)}
\bigl(\ull\t _i(t) \bigr) \bigr)\\
&&\hspace*{62pt}\qquad\quad{}-\psi' p\pmatrix{t
\cr
i-1} \bigl( \bigl(1-\psi' \bigr)|A|p \bigr)^{i-1}\bigl(1-|A|p\bigr)^{t-i+1}.
\end{eqnarray*}

One may think of the sum term to be the error generated by
approximating $\Pr(\mathcal{C}_{(x,y)}(\ull\t_i(t)))$ with the binomial
expression.
Since $\Pr(\mathcal{C}_{(x,y)}(\ull\t_i(t))) \leq(1-\Row)^i$ and
$(1-\psi') \leq(1-\Row)$, the absolute value of the sum term is at most
%
\begin{equation}
(1-\Row)^{(k \ln n)^5}{\Biggl(1+ \sum_{i=(k \ln n)^5}^{t}
\pmatrix {t
\cr
i-1}\bigl(|A|p\bigr)^{i-1}\bigl(1-|A|p\bigr)^{t-i+1}
\Biggr)}.\label{f5f5gs56dgo6p}
\end{equation}
Furthermore, $O((1-\Row)^{(k \ln n)^5})= O(n^{-\Omega(k^5 (\ln
n)^4)})$, and
\begin{eqnarray*}
\sum_{i=(k \ln n)^5}^{t}\pmatrix{t
\cr
i-1}\bigl(|A|p\bigr)^{i-1}\bigl(1-|A|p\bigr)^{t-i+1} &\leq & \sum
_{i=(k \ln n)^5}^{t}\biggl(\frac{et} {i}^i
\biggr) \bigl(|A|p\bigr)^{i}
\\
&\leq& \sum_{i=(k \ln n)^5}^{\infty}\biggl(
\frac{et|A|p} {(k
\ln n)^5}^i\biggr)
\\
&=&O{\biggl(\biggl(\frac{k^2t} {(k \ln n)^5n}^{(k \ln n)^5}\biggr)\biggr)}.
\end{eqnarray*}
Since $t \leq k^2 n \ln n$, we have $\frac{k^2t}{(k \ln n)^5n}\leq
\frac
{1}{k (\ln n)^4}$. Therefore, \eqref{f5f5gs56dgo6p} is $O(n^{-\Omega
(k^5 (\ln n)^4)})$.
Hence, when $k \leq n^\epsilon$ for small enough $\epsilon$,
\begin{eqnarray*}
\Pr \bigl(\mathcal{C}_{(x,y)}(t) \bigr)&=&\psi'p \bigl(1-
\psi'|A|p \bigr)^t - O\biggl(\frac {1} {n^{\Omega(\sqrt{n})}}\biggr)+
O \biggl(\frac{1} {n^{\Omega(k^5 (\ln n)^4)}}\biggr)
\\
&=&\psi'p \bigl(1-\psi'|A|p \bigr)^t + O
\biggl(\frac{1} {n^{\Omega(k^5 (\ln n)^4)}}\biggr)
\\
&=&{\biggl(1+O{\biggl(+\frac{(1-\psi'|A|p)^{-t}}{n^{\Omega(k^5 (\ln
n)^4)}}\biggr)}\biggr)}\psi'p
\bigl(1- \psi'|A|p \bigr)^t.
\end{eqnarray*}
Since $t \leq k^2 n \ln n$, we have for some constant $C$,
\begin{eqnarray*}
\bigl(1-\psi'|A|p \bigr)^{-t}&=&O \bigl(
\bigl(1-Ck^2/n \bigr)^{-k^2n\ln n} \bigr)=O \bigl(
\bigl(1-Ck^2/n \bigr)^{-({n}/{(Ck^2)})k^4\ln n} \bigr)\\
&=&n^{O(k^4)}.
\end{eqnarray*}
Therefore, we have $\Pr(\mathcal{C}_{(x,y)}(t))=\psi'p(1-\psi'|A|p)^t$.

Consider the correcting factor of \eqref{psi'def}; this can be absorbed
in the above correcting factors, as well as into $p$. We are finally
left with
\[
\Pr \bigl(\mathcal{C}_{(x,y)}(t) \bigr)=\psi p\bigl(1-|A|\psi
p\bigr)^t.
\]
\upqed\end{pf*}


\section{\texorpdfstring{Concluding remarks.}{Concluding remarks}}\label{ConclusionSection}
\subsection{\texorpdfstring{Comments on proof strategy.}{Comments on proof strategy}}

In Section~\ref{SIRRho1}, one may wonder why we do not employ the more
obvious strategy of simply setting the infectious period to be infinite
(i.e., letting it run as an SI process) then deleting from $\Upsi$
edges with weight exceeding the original infectious period $\xi$. We
could then have continued where Section~\ref{SIModel} left off, proving
Theorem~\ref{epidemic} without the need for the intervening material of
Section~\ref{SIRRho1}. Unfortunately, this approach fails, as can be
demonstrated by the following example: Suppose there are four particles
$a,b,c,d$ with $a$ being the initial infective. Suppose the following
(\textit{particlepair}, \textit{timestep}) meetings take place (the repetition of
$cd$ is intentional): $(ab,9)$, $(ad, 11)$, $(bc, 18)$, $(cd, 22)$,
$(cd, 27)$, $(ac, 100)$, $(bd, 100)$. The SI-based interaction graph
would have the following edge weights: $w_{\Upsi}(a,b)=9$, $w_{\Upsi
}(a,c)=100$, $w_{\Upsi}(a,d)=11$, $w_{\Upsi}(b,c)=9$, $w_{\Upsi
}(b,d)=91$, $w_{\Upsi}(c,d)=11$. Therefore, if $\xi=10$ and we had a
rule to remove any edge with weight greater than this, it would leave
$d$ isolated, suggesting it does not get infected. However, it clearly
gets infected by the chain $a \leadsto b \leadsto c \leadsto d$.

In contrast, under the current scheme, we get weights $w_{\Upsi
}(a,b)=9$,\break $w_{\Upsi}(a,c)=100$, $w_{\Upsi}(a,d)=11$, $w_{\Upsi
}(b,c)=9$, $w_{\Upsi}(b,d)=91$, $w_{\Upsi}(c,d)=4$. Deleting edges, we
see that all particles are in the connected component of $a$, meaning
they all get infected.

\subsection{\texorpdfstring{Extensions.}{Extensions}}
One obvious extension is generalising Theorem~\ref
{completiontimetheorem} to include the case $\xi< \infty$. This would
require the maximum weighted distance from $x_0$ to other vertices in
$\mathcal{C}_{f_{\xi}(\Lambda)}$. For this purpose, it may be possible
to exploit recent results such as \cite{ADL}.

Another obvious extension would be making the infectious period random,
independently for each particle. Of course, we would not be able to use
the current strategy of deleting edges that exceed a particular finite
weight, and it would appear that the techniques in this paper do not
readily extend to be able to cope with this setting. If one were to
relax the model to allow infectious periods to be associated with
particle pairs rather than particles themselves, this would correspond
to i.i.d. random cut-off thresholds on edges. However, it would be
difficult to justify an interpretation of this model. It would probably
be more fruitful to aim to calculate rough bounds on $M_k$ rather than
precise values we currently get.

Finally, one may consider other graph models. In particular, random
graphs of a prescribed degree sequence generalise random regular
graphs, so would seem an obvious extension. Such a model was studied in
\cite{CovDS}, and it would seem that some of the results and techniques
in that paper could find use in a multiple walks setting.

Another setting is $d$-dimensional grids as per \cite{Upfal} and \cite
{Lam}. Those papers study the SI model, getting results on broadcast
time, which is equivalent to our completion time $T_k$. Studying the
size of the outbreak in the SIR model is a natural avenue for
investigation. The techniques in this paper would not be amenable to
that setting, due to the very different nature of these families of
graphs. For a pertinent example, random regular graphs are mostly
locally treelike with short cycles being far from each other. Grids on
the other hand, have many short cycles. Thus local behaviours of walks
in the mixing times, (as well as the mixing times themselves) will be
quite different.


\begin{appendix}\label{app}
\section*{Appendix}

\begin{pf*}{Proof of Lemma~\ref{TimeDistPsi}}
{(i)${}\Rightarrow{}$(ii)}
It is straightforward to show that for a particle $y$, the infection
time $t(y)$ is the weighted distance in $f_\xi(\Upsi')$ between $x_0$
and~$y$, which will be $\infty$ if there is no path. Therefore,
{(i)} implies that for an edge $e=(x,y)$, $t(e)=d_{f_\xi(\mathbf
{z})}(e)$, and by the construction of $\Upsi'$, {(ii)} follows.

{(ii)${}\Rightarrow{}$(i)}
Order particles in $\mathcal{C}_{f_\xi(\mathbf{z})}$ by their weighted
distance from $x_0$ in $f_\xi(\mathbf{z})$: $x_0=x_{(0)}, x_{(1)},
\ldots, x_{(r-1)}$ where $r=|\mathcal{C}_{f_\xi(\mathbf{z})}|$ and $i<j
\Rightarrow d_{f_{\xi}(\mathbf{z})}(x_{(i)}) \leq d_{f_{\xi}(\mathbf
{z})}(x_{(j)})$. We shall prove that $t(x_{(i)})=d_{f_{\xi}(\mathbf
{z})}(x_{(i)})$ for $0\leq i\leq r-1$.

Clearly this proposition holds for $x_{(0)}$. Suppose for all $i \leq
N-1$, $t(x_{(i)})=d_{f_{\xi}(\mathbf{z})}(x_{(i)})$. If $N=r$, we are
done. Otherwise $N<r$. Let $x_{(M)}$ be a neighbour of $x_{(N)}$ on a
shortest path from $x_0$ to $x_{(N)}$. Since $d_{f_{\xi}(\mathbf
{z})}(x_{(M)}) < d_{f_{\xi}(\mathbf{z})}(x_{(N)})$, $M < N$, so by the
induction hypothesis, $t(x_{(M)})=d_{f_{\xi}(\mathbf
{z})}(x_{(M)})=d_{f_{\xi}(\mathbf{z})}(e)$ where $e=(x_{(M)},
x_{(N)})$. By {(ii)} this implies $t(x_{(N)}) \leq d_{f_{\xi
}(\mathbf{z})}(e)+w_\mathbf{z}(e)=d_{f_{\xi}(\mathbf{z})}(x_{(N)})$.\vspace*{1pt}

Now consider the chain of infections starting at $x_0$ that led to
$x_{(N)}$ being infected. Let $x_{(j)}$ be the first in the chain where
$j>N-1$, and suppose it got infected by $x_{(i)}$, $i \leq N-1$. Then
by {(ii)}, $t(x_{(j)}) = d_{f_{\xi}(\mathbf
{z})}(x_{(i)})+w_{\mathbf{z}}(x_{(i)}, x_{(j)})\geq d_{f_{\xi
}(\mathbf
{z})}(x_{(j)}) \geq d_{f_{\xi}(\mathbf{z})}(x_{(N)})$, implying
$t(x_{(N)}) \geq d_{f_{\xi}(\mathbf{z})}(x_{(N)})$. Hence
$t(x_{(N)})=d_{f_{\xi}(\mathbf{z})}(x_{(N)})$ for $0\leq i\leq r-1$.
\end{pf*}

\begin{pf*}{Proof of Lemma~\ref{probConv}}
Let $A(\mathbf{z})=V(\mathcal{C}_{f_\xi(\mathbf{z})})$, $B(\mathbf
{z})=\mathcal{P}\setminus A(\mathbf{z})$ and $E(\mathbf
{z})=E(A(\mathbf
{z}))\cup E(A(\mathbf{z})\dvtx B(\mathbf{z}))$. Let $E_{\Lambda}(\mathbf
{z})=E_{\Lambda}(A(\mathbf{z}))\cup E_{\Lambda}(A(\mathbf
{z})\dvtx B(\mathbf
{z}))$, and similarly for $\Upsi'$. Since $\Lambda$ is good w.h.p., by
Lemma~\ref{ineractionapproxPsi},
\begin{eqnarray*}
1-o(1)&=&\sum_{E(\mathbf{z})\dvtx\mathbf{z}\ \mathrm{good}}\Pr \bigl(E_{\Lambda
}(
\mathbf{z})
=E(\mathbf{z}) \bigr)\\
&=& \bigl(1-o(1) \bigr)\sum
_{E(\mathbf{z})\dvtx \mathbf{z}\ \mathrm{good}}\Pr \bigl(E_{\Upsi'}(\mathbf{z})=E(\mathbf{z})
\bigr).
\end{eqnarray*}
Observe $\mathcal{C}_{f_\xi(\mathbf{z})} =\mathcal{C}_{f_\xi
(E(\mathbf
{z}))}$. Let $\mathbf{1}_{\mathcal{E}_n}(\mathcal{C}_{f_\xi
(E(\mathbf
{z}))})$ be the indicator for $\mathcal{C}_{f_\xi(E(\mathbf{z}))}\in
\mathcal{E}_n$.
\begin{eqnarray*}
\Pr(\mathcal{C}_{f_\xi(\Lambda)} \in\mathcal{E}_n) &=& \sum
_{E(\mathbf{z})}\Pr \bigl(E_{\Lambda}(\mathbf{z})=E(
\mathbf {z}) \bigr)\mathbf{1}_{\mathcal{E}_n}(\mathcal{C}_{f_\xi(E(\mathbf
{z}))})
\\
&=& o(1)+ \sum_{E(\mathbf{z})\dvtx \mathbf{z}\ \mathrm{good}}\Pr \bigl(E_{\Lambda
}(
\mathbf{z})=E(\mathbf{z}) \bigr)\mathbf{1}_{\mathcal{E}_n}(\mathcal
{C}_{f_\xi(E(\mathbf{z}))})
\\
&=& o(1)+ \bigl(1-o(1) \bigr)\sum_{E(\mathbf{z})\dvtx \mathbf{z}\ \mathrm{good}}\Pr
\bigl(E_{\Upsi'}(\mathbf{z})=E(\mathbf{z}) \bigr)\mathbf{1}_{\mathcal
{E}_n}(
\mathcal{C}_{f_\xi(E(\mathbf{z}))})
\\
&=& o(1)+ \bigl(1-o(1) \bigr)\sum_{E(\mathbf{z})}\Pr
\bigl(E_{\Upsi'}(\mathbf {z})=E(\mathbf{z}) \bigr)\mathbf{1}_{\mathcal{E}_n}(
\mathcal{C}_{f_\xi
(E(\mathbf{z}))})
\\
&=& o(1)+ \bigl(1-o(1) \bigr)\Pr(\mathcal{C}_{f_\xi(\Upsi')} \in
\mathcal{E}_n).
\end{eqnarray*}
Hence if $\Pr(\mathcal{C}_{f_\xi(\Lambda)} \in\mathcal
{E}_n)\rightarrow p_c$ as $n \rightarrow\infty$, then $\Pr(\mathcal
{C}_{f_\xi(\Upsi')} \in\mathcal{E}_n)\rightarrow p_c$ as $n
\rightarrow\infty$. Corollary~\ref{UpsiesAREequal}({i}) says
$\Upsi=\Upsi'$ w.h.p., and so $\Pr(\mathcal{C}_{f_\xi(\Upsi)} \in
\mathcal
{E}_n)\rightarrow p_c$ as $n \rightarrow\infty$
\end{pf*}

\begin{pf*}{Proof of Lemma~\ref{GoodLemExt}}
{(a)} In the proof of Lemma~\ref{UpsiesEqual} we derived
probabilities for each of these two types of ``failures'': {(i)}
more than one pair of active particles being within distance $d=\alpha
(\ln\ln n +\ln k)$ of each other when some pair meet after the
(extended) mixing time, and {(ii)} any pair meeting within
$\ell
$ steps\vspace*{1pt} when they all start with distance at least $d$ from each other.
These had probabilities $O(1/n^{\Omega(1)})$ and $O(1/(k \ln
n)^{\Omega
^+(1)})$, respectively. Each type of failure was considered over at
most ${k\choose2}$ particle pairs over at most ${k\choose2}$ meetings
(interactions), so taking the union bound over all of them, the above
probabilities were still $o(1)$.

Now consider the $\Row<1$ setting. After the extended mixing time we
wait for a meeting of an active pair. Say the active pair that meets is
$(x,y)$. We allow $(x,y)$ to have $T$ steps to interact, before letting
the system mix again for $2(T+T^3)$ steps, after which we wait for
another active pair meeting. Thus, if there are no failures as
described above, these first meetings after the extended mixing times
are $\ell$-distinct. Let the random variable $X$ count the total number
of such meetings over the course of the process, over all active
particle pairs. The probability of an interaction in the $T$-length
window is at least $\Row$. Since $\Row$ is constant, $\E
[X]=O\bigl({k\choose
2}/\Row\bigr)=O(k^2)$. Hence, the expected number of failures over all
particle pairs over all $X$ meetings is bounded by $\E[X(1/n^{\Omega
(1)}+1/(k \ln n)^{\Omega^+(1)})]=o(1)$.

{(b)}
We can assume $\t=0$. Let $E_1$ be the event that $x$ and $y$ interact
in the period $[0,T-1]$, and let $E_2$ be the event that $x$ and $y$
meet in the period $[T, \ell]$. As per Lemma~\ref{2-interact}, $\Pr
(E_1)=\psi'$. Then $\Pr(E_2 | \overline{E}_1) \leq\Pr(E_2)/\Pr
(\overline{E}_1)=\Pr(E_2)/(1-\psi')$.

Thus, since $\psi'$ is almost a constant, $\Pr(E_2| \overline{E}_1)$ is
of the same order as $\Pr(E_2)$. The rest of the proof is similar to
part {(a)}.
\end{pf*}
\end{appendix}




\printaddresses

\end{document}